\documentclass[11pt]{amsart}
\usepackage{latexsym,amssymb,amsfonts,amsmath,graphicx,fullpage,url,color}
\usepackage[utf8]{inputenc}
\usepackage[vcentermath,enableskew]{youngtab}
\usepackage[all]{xy}
\usepackage{verbatim}
\usepackage{ytableau}
\usepackage{makecell}
\usepackage[export]{adjustbox}
\usepackage{multicol}
\usepackage{xcolor}

\usepackage{mathdots,verbatim}
\usepackage{tikz, calc, ifthen}
\usetikzlibrary{calc, shapes, backgrounds,arrows,positioning}
\tikzset{>=stealth',
  head/.style = {fill = white, text=black}, 
  pil/.style={->,thick},
  junct/.style = {draw,circle,inner sep=0.5pt,outer sep=0pt, fill=black}
  }
\definecolor{light-gray}{gray}{0.85}
\definecolor{dark-gray}{gray}{0.25}

\usepackage{tikz,tkz-graph}

\newcounter{x}
\newcounter{y}
\newcounter{z}

\newcommand\xaxis{210}
\newcommand\yaxis{-30}
\newcommand\zaxis{90}

\newcommand\topside[3]{
  \fill[fill=white, draw=black,shift={(\xaxis:#1)},shift={(\yaxis:#2)},
  shift={(\zaxis:#3)}] (0,0) -- (30:1) -- (0,1) --(150:1)--(0,0);
}

\newcommand\leftside[3]{
  \fill[fill=black, draw=black,shift={(\xaxis:#1)},shift={(\yaxis:#2)},
  shift={(\zaxis:#3)}] (0,0) -- (0,-1) -- (210:1) --(150:1)--(0,0);
}

\newcommand\rightside[3]{
  \fill[fill=gray, draw=black,shift={(\xaxis:#1)},shift={(\yaxis:#2)},
  shift={(\zaxis:#3)}] (0,0) -- (30:1) -- (-30:1) --(0,-1)--(0,0);
}

\newcommand\cube[3]{
  \topside{#1}{#2}{#3} \leftside{#1}{#2}{#3} \rightside{#1}{#2}{#3}
}

\newcommand\planepartition[1]{
 \setcounter{x}{-1}
  \foreach \a in {#1} {
    \addtocounter{x}{1}
    \setcounter{y}{-1}
    \foreach \b in \a {
      \addtocounter{y}{1}
      \setcounter{z}{-1}
      \foreach \c in {0,...,\b} {
        \addtocounter{z}{1}
      \ifthenelse{\c=0}{\setcounter{z}{-1},\addtocounter{y}{0}}{
        \cube{\value{x}}{\value{y}}{\value{z}}
      }
    }
  }
 }
}

\usetikzlibrary{arrows}
\usetikzlibrary{shapes}

\newtheorem{theorem}{Theorem}[section]
\newtheorem{proposition}[theorem]{Proposition}
\newtheorem{lemma}[theorem]{Lemma}
\newtheorem{corollary}[theorem]{Corollary}
\newtheorem{conjecture}[theorem]{Conjecture}
\newtheorem*{thmbdry}{Theorem~\ref{thm:bdry_ones}}
\newtheorem*{thmmagogbij}{Theorem~\ref{prop:magog_bij}}
\newtheorem*{thm132}{Theorem~\ref{thm:132}}
\newtheorem*{thmmagogineq}{Theorem~\ref{thm:magogineq}}
\newtheorem*{thminv}{Theorem~\ref{thm:inv}}
\newtheorem*{thmbtpineq}{Theorem~\ref{thm:btpineq}}
\newtheorem*{thmbtpfacets}{Theorem~\ref{thm:btpfacets}}
\theoremstyle{definition}
\newtheorem{definition}[theorem]{Definition}
\newtheorem{example}[theorem]{Example}
\newtheorem{remark}[theorem]{Remark}
\newtheorem*{magogdef}{Definition~\ref{ineq:m}}

\newcommand{\ZZ}{\mathbb{Z}}
\newcommand{\RR}{\mathbb{R}}

\newcommand{\ASM}{\mathrm{ASM}}
\newcommand{\TSSCPP}{\mathrm{TSSCPP}}
\newcommand{\magog}{\mathrm{Magog}}
\newcommand{\sign}{\mathrm{Sign}}
\newcommand{\inv}{\mathrm{inv}}
\newcommand{\posinv}{\mathrm{posinv}}
\newcommand{\BTP}{\mathrm{BTP}}

\numberwithin{equation}{section}

\title{Totally symmetric self-complementary plane partition matrices and related polytopes}
\author{Vincent Holmlund}
\author{Jessica Striker} 
\address{Department of Mathematics, North Dakota State University}
\email{vincent.holmlund@ndsu.edu}
\email{jessica.striker@ndsu.edu}
\date{\today}

\begin{document}

\begin{abstract}
Plane partitions in the totally symmetric self-complementary symmetry class (TSSCPPs) are known to be equinumerous with $n\times n$ alternating sign matrices, but no explicit bijection is known. In this paper, we give a bijection from these plane partitions to $\{0,1,-1\}$-matrices we call magog matrices, some of which are alternating sign matrices. We explore enumerative properties of these matrices related to natural statistics such as inversion number and number of negative ones. We then investigate the polytope defined as their convex hull. We show that all the magog matrices are extreme and give a partial inequality description. Finally, we define another TSSCPP polytope as the convex hull of TSSCPP boolean triangles and determine its dimension, inequalities, vertices, and facets.
\end{abstract}

\maketitle

\section{Introduction}
Alternating sign matrices (ASMs) are certain $\{0,1,-1\}$-matrices that form an interesting superset of permutations. They have  nice properties from many perspectives. 
Enumerative results include a lovely counting formula (conjectured in~\cite{MILLS1983340} and proved in 
\cite{zeilberger1994proof,kuperberg1997proof}) and  refined enumerations~\cite{behrend.multiply.refined}. 
Geometrically, their convex hull gives an intriguing polytope~\cite{Behrend2007HigherSA,Striker2007TheAS} with many  properties mirroring those of the Birkhoff polytope. Poset-theoretically, they form a distributive lattice whose subposet of join irreducibles is a tetrahedron~\cite{ELKP92,Striker.tetrahedralposet}. 
Algebraic results include a theorem that the ASM lattice is the MacNeille completion of the strong Bruhat order on permutations~\cite{Lascoux.Schutzenberger} and a very recent result in which ASMs appear in $SL_4$ web basis equivalence classes~\cite{Gaetz.Pechenik.Pfannerer.Striker.Swanson:4row}.

The most mysterious property of ASMs is that they are equinumerous with other sets of objects with whom no explicit bijection is known. In this paper, we transform one such set of objects, totally symmetric self-complementary plane partitions (TSSCPPs), into matrices with entries in $\{0,1,-1\}$ we call \emph{magog matrices}. We then study these from enumerative and geometric perspectives, comparing with analogous results on ASMs.

Before describing magog matrices, we first mention a superset of both ASMs and magog matrices: the \emph{sign matrices} of Aval~\cite{aval2007keys}. Definition~\ref{def:sq_sign} of the present paper isolates a subset of sign matrices that contains both ASMs and magog matrices. In Section~\ref{sec:matrices}, we define magog matrices as square sign matrices with one additional set of inequalities.

\begin{magogdef}
    An $n\times n$ \emph{magog matrix} $A=(a_{ij})$ is a square sign matrix such that:
    \begin{align*}          
     \sum_{j'=1}^{j}a_{i+1,j'}+ \sum_{i'=1}^{i+1}a_{i',j+1}-\sum_{i'=1}^{i}a_{i'j} \geq 0 & \text{ for all }
     1 \leq i \leq n-2, 
     1 \leq j \leq n-2.  
     \end{align*}
\end{magogdef}

Our first main result is the theorem below, which shows a bijection between magog matrices and TSSCPPs.

\begin{thmmagogbij}
    The set of $n\times n$ magog matrices is in explicit bijection with the set of magog triangles of order $n$, and therefore with the set of all TSSCPPs in a $2n\times 2n\times 2n$ box.
\end{thmmagogbij}

The following enumerative results are the other main theorems of Section~\ref{sec:matrices}.

\begin{thm132}
    The magog matrices of order $n$ with no negative ones are the $132$-avoiding permutation matrices.
\end{thm132}

\begin{thmbdry} 
Let $\magog_{n}(i,j)$ be the set of magog matrices of order $n$ with a one in row $i$ column $j$. The following enumerative identities hold:
    \begin{enumerate}
        \item For all $n\geq 1$,  $|\magog_{n}(1,1)|=1$,
        \item   For all $n>1$, $|\magog_{n}(n,1)|=|\magog_{n}(n,2)|=|\magog_{n-1}|=\displaystyle\prod_{j=0}^{n-2}\frac{(3n-2)!}{(n-1+j)!}$,
        \item For all $n>1$, $|\magog_{n}(1,n)|=|\magog_{n}(1,n-1)|$,
        \item For $n\geq 1$, $|\magog_{n}(2,1)|=\displaystyle\frac{1}{n+1}\binom{2n}{n}-1$, and
        \item For all $n\geq 1$,
    $|\magog_{n}(1,2)|=2^{n-1}-1$.
    \end{enumerate}
\end{thmbdry}

\begin{thminv}
Let $\magog_n(k~\inv)$ (resp. $\magog_n(k~\posinv)$) denote the set of $n\times n$ magog matrices with inversion number (resp. positive inversion number) $k$. 
The following enumerations hold:
\begin{enumerate}
    \item  For all $n\geq 1$, $|\magog_n(0~\inv)|=|\magog_n(0~\posinv)|=1$,
    \item For all $n\geq 2$, $|\magog_n(1~\inv)|=1$, 
    \item For all $n\geq 2$, $|\magog_n(\binom{n}{2}-1~\posinv)|=n-1$,
    \item For all $n\geq 2$, $|\magog_n(\binom{n}{2}~\inv)|=|\magog_n(\binom{n}{2}~\posinv)|=1$, and
    \item For all $n\geq 3$, $|\magog_n(2~\inv)|=n+1$.
\end{enumerate}
\end{thminv}

In Section~\ref{sec:magogpoly}, we take the convex hull of $n\times n$ magog matrices to form a polytope we call $\TSSCPP(n)$. We prove a partial inequality description in the theorem listed below and also show that each magog matrix is extreme in Theorem~\ref{thm:magogextreme}.

    \begin{thmmagogineq}
       The following inequalities hold for any matrix $A=(a_{ij}) \in \TSSCPP(n)$: 
        \begin{align*}      
        \sum_{j'=1}^{j}a_{i+1,j'}+\sum_{i'=1}^{i+1}a_{i',j+1} \geq 1
       &\text{ for all }
        1 \leq i \leq n-2,
        1 \leq j \leq n-2,
        i+j \geq n-1,  \\ 
        \sum_{j'=1}^{j+1}a_{2,j'}+\sum_{j'=j+1}^{n}a_{1,j'} \geq 1  
        &\text{ for all } 
        1 \leq j \leq n-3, \\
\sum_{i'=1}^{i+1}a_{i',2}+\sum_{i'=i+1}^{n}a_{i',1} \geq 1        
        &\text{ for all }
        1 \leq i \leq n-3.    
         \end{align*}    
    \end{thmmagogineq}

We end the paper 
with a different geometric perspective on TSSCPPs, by forming a polytope from alternative data: the associated boolean triangles. This polytope turns out to be better behaved; we prove a complete inequality description and count the facets in the theorems listed below. We also prove each boolean triangle is extreme in Theorem~\ref{thm:btpextreme}.

    \begin{thmbtpineq}
     The TSSCPP boolean triangle polytope $\BTP(n)$ is equal to the set of all triangular arrays $(b_{i,n-j})$ for $1 \leq j \leq i \leq n-1$ such that: 
        \begin{align*}
        b_{i,n-j} &\geq 0, \\
        b_{i,n-j} &\leq 1, \\
    1+\sum_{k=j+1}^{i}b_{k,n-j-1} &\geq \sum_{k=j}^{i}b_{k,n-j}    \text{ for all } 1 \leq j < i \leq n-1. 
    \end{align*}
    \end{thmbtpineq}

\begin{thmbtpfacets}
        The number of facets of $\BTP(n)$ is $\displaystyle\frac{(n-1)(3n-2)}{2}$.
\end{thmbtpfacets}

This paper is organized as follows. In Section~\ref{sec:back}, we give background on ASMs, sign matrices, and related work on their polytopes. We then discuss background on TSSCPPs, including their characterization as magog triangles and boolean triangles.
In Section~\ref{sec:matrices}, we define magog matrices, show they are in bijection with TSSCPPs, and study their enumerative properties with respect to number of negative ones, boundary ones in specified entries, and inversion numbers. In Section~\ref{sec:magogpoly}, we define and study the polytope formed as the convex hull of $n\times n$ magog matrices. Section~\ref{sec:btp} similarly studies the polytope formed as the convex hull of TSSCPP boolean triangles of order $n$. Many subsections include data tables and a closing remark comparing its results with similar results on alternating sign matrices.

\section{Background}
\label{sec:back}
In this section, we give needed background on polytopes and TSSCPPs, as well as ASMs and their related polytope. Readers familiar with these topics may still want to consult Subsection~\ref{subsec:sign}, which discusses the less familiar \emph{sign matrices} and their polytope.
\subsection{Polytopes}
\label{subsec:backpoly}
We begin with the two ways to define polytopes and their equivalence, as discussed in \cite{ziegler2012lectures}.
    A \emph{$\mathcal{V}$-polytope} $P \subseteq \RR^n$ is a convex hull of a finite set of points in $\RR^n$.
    An \emph{$\mathcal{H}$-polytope} $P \subseteq \RR^n$ is the bounded intersection of a finite number of half spaces in $\RR^n$.
It is an important theorem in convex geometry that
    any $\mathcal{H}$-polytope is also a $\mathcal{V}$-polytope and vice versa.

The Birkhoff polytope, originally studied by  Birkhoff and von Neumann~\cite{BIRKHOFF1946,VonNeumann1953}, is the convex hull of all $n\times n$ permutation matrices. Its description as an $\mathcal{H}$-polytope is as the set of all $n\times n$ matrices with non-negative entries whose row and column sums equal one. 
The matrix polytopes discussed in this paper all relate to the Birkhoff polytope.

The \emph{$f$-vector} of a polytope is a vector whose $i$-th coordinate is the number of faces of dimension $i-1$. A \emph{facet} of a polytope is a face of dimension one less than the polytope itself.

\subsection{Alternating sign matrices and the ASM polytope} 
\label{subsec:backASM}
In this and the next subsections, we review definitions of certain matrices and prior work on their polytopes.
\begin{definition}[\protect{\cite[Definition 1]{MILLS1983340}}]
    An \emph{alternating sign matrix (ASM)} is a square $\{0,1,-1\}$-matrix where the rows and columns sum to one and the nonzero entries in each row and column alternate in sign. Let $\ASM_n$ denote the set of $n\times n$ ASMs.
\end{definition} 

The left seven matrices in Figure~\ref{fig:8sq} are the matrices of $\ASM_3$.
\smallskip

Alternating sign matrices were initially studied by Robbins and Rumsey in relation to the lambda determinant \cite{bressoud_1999}. The enumeration formula for $n\times n$ ASMs,
\begin{equation}
\label{asmform}
\prod_{j=0}^{n-1}\frac{(3n+1)!}{(n+j)!},
\end{equation}
was initially conjectured by Mills, Robbins, and Rumsey \cite{MILLS1983340} and was subsequently proven independently by Zeilberger \cite{zeilberger1994proof}, Kuperberg \cite{kuperberg1997proof}, and  Fischer \cite{FISCHER2007253}. For more details on the history of ASMs, see \cite{bressoud_1999}. 

Striker studied ASMs geometrically by taking the convex hull of $n \times n$ ASMs to define  the \emph{alternating sign matrix polytope}~\cite{Striker2007TheAS}. She proved analogous theorems to those known on the Birkhoff polytope on topics such as the facets, vertices, inequality description, and dimension. This polytope was independently defined and studied by Behrend and Knight. They studied the vertices and inequality description. In addition, they studied lattice points in the $r$th-dilate of the ASM polytope, termed \emph{higher spin ASMs}~\cite{Behrend2007HigherSA}.

\subsection{Sign matrices and sign matrix polytopes}
\label{subsec:sign}
One generalization of ASMs is the set of \emph{sign matrices}, defined and shown to be in bijection with semistandard Young tableaux by Aval~\cite{aval2007keys}.
Sign matrices are a superset of ASMs in that they are not necessarily square, the row and column sums need not be one, and the row partial sums must be non-negative as opposed to the stronger ASM condition of being between zero and one. 

We will use a special case of sign matrices in this paper that we call \emph{square sign matrices}. These are sign matrices whose shape is a square and also whose rows each sum to one.
\begin{definition}
\label{def:sq_sign}
        An $n\times n$ \emph{square sign matrix} is a square $\{0,1,-1\}$-matrix $A=(a_{ij})$ such that:
    \begin{align}        
    \sum_{i'=1}^{n}a_{i'j}=1 & \text{ for all } 1 \leq j \leq n, \label{eq:1magog}\\         
    \sum_{j'=1}^{n}a_{ij'}=1 & \text{ for all } 1 \leq i \leq n, \label{eq:2magog}\\        
    0 \leq \sum_{i'=1}^{i}a_{i'j} \leq 1 & \text{ for all } 1 \leq i \leq n, 1 \leq j \leq n, \label{eq:3magog}  \\  
     \sum_{j'=1}^{j}a_{ij'} \geq 0    & \text{ for all } 1 \leq i \leq n, 1 \leq j \leq n.\label{eq:4magog}  
     \end{align}
    We call the set of all square sign matrices $\sign_n$.
\end{definition}

Alternating sign matrices are square sign matrices that additionally satisfy     
\begin{equation}
    \label{eq:ASMineq}
\sum_{j'=1}^{j}a_{ij'} \leq 1     
\text{ for all }  1 \leq i \leq n, 1 \leq j \leq n.
 \end{equation}
 See Figure~\ref{fig:8sq} for the eight matrices in $\sign_3$.
\begin{figure}[htbp]
\scalebox{.95}{
$\begin{pmatrix} 1 & 0 & 0 \\ 0 & 1 & 0\\ 0 & 0 & 1 \end{pmatrix}
 \begin{pmatrix} 0 & 1 & 0 \\   1 & 0 & 0\\ 0 & 0 & 1 \end{pmatrix}
 \begin{pmatrix} 0 & 1 & 0 \\ 0 & 0 & 1 \\ 1 & 0 & 0 \end{pmatrix}
 \begin{pmatrix} 0 & 1 & 0 \\ 1 & -1 & 1\\ 0 & 1 & 0\end{pmatrix}
 \begin{pmatrix} 0 & 0 & 1 \\ 1 & 0 & 0 \\ 0 & 1 & 0 \end{pmatrix}
 \begin{pmatrix} 0 & 0 & 1 \\ 0 & 1 & 0 \\ 1 & 0 & 0 \end{pmatrix}
 \begin{pmatrix} 1 & 0 & 0 \\ 0 & 0 & 1 \\ 0 & 1 & 0 \end{pmatrix}
 \begin{pmatrix} 0 & 0 & 1 \\ 1 & 1 & -1 \\ 0 & 0 & 1 \end{pmatrix}$}
\caption{The eight $3 \times 3$ square sign matrices. All except the last are ASMs. All except the second to last are magog matrices.}
\label{fig:8sq}
\end{figure}

The enumeration below follows from Aval's bijection to tableaux. We give a short proof.
\begin{proposition}
    The number of $n\times n$ square sign matrices is $2^{\binom{n}{2}}$.
\end{proposition}
\begin{proof}
By \cite{aval2007keys}, the set $\sign_n$ of square sign matrices is in bijection with semistandard Young tableaux of shape $(n,n-1,\ldots,2,1)$ and entries at most $n$. The enumeration then follows from Stanley's hook-content formula~\cite[Theorem 15.3]{STANLEY_PP_II}. 
\end{proof}

Solhjem and Striker studied polytopes formed as the convex hull of sign matrices for given size and row sums \cite{SOLHJEM201984}. 
They investigated the inequality descriptions, vertices, facets, and face lattice. Another generalization of ASMs is the set of \emph{partial alternating sign matrices}, sign matrices whose row sums are between zero and one. Heuer and Striker studied polytopes formed by the set of $m \times n$ partial alternating sign matrices~\cite{heuer2022partial}. They studied the inequality description, facets, and face lattice.

\subsection{TSSCPPs as magog triangles and boolean triangles}
\label{subsec:backTSSCPP}
\begin{figure}[htbp]
\begin{center}
    \includegraphics[scale=0.95]{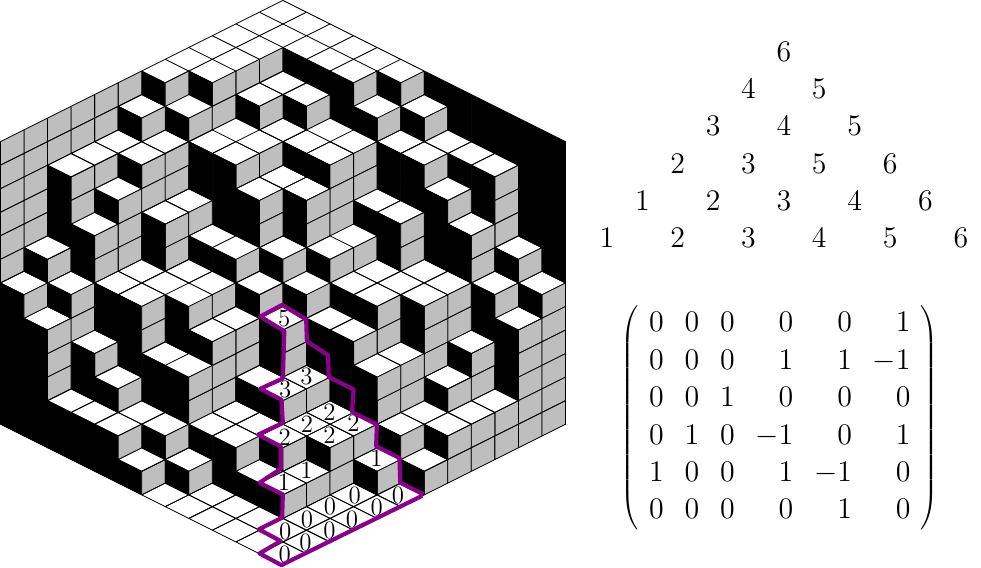}
\end{center}
\caption{Left: A TSSCPP with its fundamental domain outlined and the number of boxes in each stack listed. Top right: Its corresponding magog triangle. Bottom right: the main new object of this paper, its magog matrix. The magog triangle is obtained from the numbers in the fundamental domain by adding $i$ to column $i$ and then rotating clockwise slightly. The magog matrix construction is described in Theorem~\ref{prop:magog_bij}.}
\label{fig:bigFD}
\end{figure}

Plane partitions are three dimensional analogues of ordinary partitions. These objects were first studied by MacMahon \cite{macmahon2001combinatory}. The focus of our study is a special kind of plane partition.
\begin{definition}
    A \emph{plane partition} $\pi$ is a  set of  lattice points with positive coordinates $(i,j,k)$ such that if $(i,j,k) \in \pi$ and $1 \leq i' \leq i, 1 \leq j' \leq j, 1\leq k' \leq k$ then $(i',j',k') \in \pi$. A plane partition is \emph{symmetric} if whenever $(i,j,k) \in \pi$ then $(j,i,k) \in \pi$. A plane partition is \emph{totally symmetric} if whenever $(i,j,k) \in \pi$ all six permutations of $(i,j,k)$ are in $\pi$.
    A plane partition is \emph{self-complementary} inside a bounding box $a \times b \times c$ if it is equal to its complement in the box, that is, the collection of empty cubes in the box is the same shape as the collection of cubes in $\pi$.
    A \emph{totally symmetric self-complementary plane partition (TSSCPP)} inside a $2n \times 2n \times 2n$ box is a plane partition that is both totally symmetric and self-complementary.
\end{definition}
See Figure~\ref{fig:bigFD}, left, for an example of a TSSCPP.

 Stanley defined ten symmetry classes of plane partitions, including TSSCPPs~\cite{StanleyPP}, and proved or conjectured product formulas for all of their enumerations. Andrews proved this conjectured TSSCPP enumeration in 
  1991~\cite{ANDREWS199428}, showing that the enumeration formula for $n\times n$ ASMs \eqref{asmform} also counts $2n \times 2n \times 2n$ TSSCPPs. Zeilberger's proof of the ASM enumeration a year later established that these two sets of objects were equinumerous \cite{zeilberger1994proof}. However, an explicit bijection between the two sets has yet to be found. Some partial bijections are known~\cite{Ayyer312, Bettinelli, Biane_Cheballah_1,striker2018permutation}. The most recent one was discovered by Huang and Striker~\cite{huang2023pipe} by interpreting TSSCPPs as certain \emph{pipe dreams}.  
 
 In this paper, we will use several objects which are known to be in bijection with TSSCPPs. We now define one of them.

\begin{definition}
\label{def:magogtri}
    A \emph{magog triangle} of order $n$ is a triangular array of integers $\alpha_{ij}$ for $1 \leq i \leq n$, $n-i \leq j \leq n-1$, such that for all $i,j$: 
    \begin{align}
        \label{eq:tri_len}
        \alpha_{ij} &\leq n,\\ 
        \label{eq:tri_rowinc}
        \alpha_{ij} &< \alpha_{i,j+1},\\
        \label{eq:tri_diagplus1}
        \alpha_{ij}+1 &\geq \alpha_{i+1,j}, \text{ and}\\
        \label{eq:tribottomrow}
        \alpha_{n,j}&=j.
        \end{align}
\end{definition}

Note that \eqref{eq:tri_rowinc} and \eqref{eq:tri_diagplus1} imply:
\begin{equation}
\label{eq:tri_implied}
    \alpha_{ij} \leq \alpha_{i-1,j+1}.
\end{equation}

\begin{example} 
\label{ex:magogtriex} Below is an example of a generic magog triangle. The \rotatebox[origin=c]{35}{\textcolor{blue}{$\leq$}} symbol indicates the \eqref{eq:tri_implied} inequality, \textcolor{blue}{$<$} indicates \eqref{eq:tri_rowinc}, and \rotatebox[origin=c]{-32}{\textcolor{blue}{$+1\geq$}} indicates \eqref{eq:tri_diagplus1}. The last row entries are given by \eqref{eq:tribottomrow}.
    \[\begin{array}{ccccccccccc}
        &&  &  &  &   \alpha_{1,n-1} &  &    &  \\
         &&  &  & \alpha_{2,n-2} &  & \alpha_{2,n-1} &  &    \\ 
         &&  & \alpha_{3,n-3} &  & \alpha_{3,n-2} &  & \alpha_{3,n-1} &    \\
         && \rotatebox[origin=c]{35}{\textcolor{blue}{$\leq$}}&  &  &  \vdots  &&  &   \rotatebox[origin=c]{-32}{\textcolor{blue}{$+1\geq$}}  \\
          &\alpha_{n-1,1} &  & \alpha_{n-1,2} &   \ldots &\textcolor{blue}{<} & \ldots & \alpha_{n-1,n-2} & & \alpha_{n-1,n-1} \\
          1 && 2 && & \ldots &&& n-1 && n
     \end{array}\]
\end{example}

See Figure~\ref{fig:bigFD} for an example of the magog triangle associated to a given TSSCPP. 
\smallskip

Mills, Robbins, and Rumsey showed that triangles equivalent to these were in bijection with TSSCPPs \cite{MILLS1986277}, as they relate to the number of cubes stacked in each position of the fundamental domain. Zeilberger later gave the name magog triangle to another form equivalent to these same triangles in his proof of the ASM enumeration conjecture \cite{zeilberger1994proof}. The definition given above will be most convenient for our purposes and is used in \cite{Striker.tetrahedralposet,striker2018permutation,huang2023pipe}. See \cite[Section 2.2]{striker2018permutation} for details on the bijection and Figure~\ref{fig:bigFD} for an example.

\begin{theorem}[\protect{\cite[Theorem 1]{MILLS1986277}}, see also \protect{\cite[Proposition 2.10]{striker2018permutation}}]
\label{thm:magogtri}
    The magog triangles of order $n$ are in bijection with the set of all  TSSCPPs in a $2n \times 2n \times 2n$ box. 
\end{theorem}

The next object in bijection with TSSCPPs will be used in Section~\ref{sec:btp}, but we will see inequalities similar to these throughout the paper. These triangles represent the lattice paths that bound each layer of the fundamental domain; see \cite[Section 2.2]{striker2018permutation} for details and history.
\begin{definition}[\protect{\cite[Definition 2.12]{striker2018permutation}}]
\label{def:bool}
      A \emph{TSSCPP boolean triangle} of order $n$ is a triangular array 
      $b_{i,n-j}$ for $1 \leq j \leq i \leq n-1$ with entries in
$\{0, 1\}$ such that the partial sums satisfy the following inequality for all $1 \leq j < i \leq n - 1$:
    \begin{align}   
    \label{eq:booltri}
    1+\sum_{k=j+1}^{i}b_{k,n-j-1} \geq \sum_{k=j}^{i}b_{k,n-j} 
    \end{align}
    We call this last inequality  the $(i,j)$-inequality.
\end{definition}
\begin{example} 
\label{ex:bool1} Below is an example of the indexing in a TSSCPP boolean triangle.
    \[\begin{array}{ccccccccc}
        &  &  &  &   b_{1,n-1} &  &    &  \\
         &  &  & b_{2,n-2} &  & b_{2,n-1} &  &    \\ 
         &  & b_{3,n-3} &  & b_{3,n-2} &  & b_{3,n-1} &    \\
         & &  &  &  \vdots  &&  &    \\
          b_{n-1,1} &  & b_{n-1,2} & &  \ldots &  & b_{n-1,n-2} & & b_{n-1,n-1} \\
     \end{array}\]
\end{example}
\begin{example}
\label{ex:bool2}
    Below, left, is a nonexample of a TSSCPP boolean triangle. The $(3,1)$-inequality is not satisfied, since $\sum_{k=1}^{3}b_{k,n-1}=3$ while $\sum_{k=2}^{3}b_{k,n-2}=1$.
    The right triangle is the TSSCPP boolean triangle corresponding to the TSSCPP of Figure~\ref{fig:bigFD}; note all $(i,j)$-inequalities are satisfied. 
\[
     \begin{array}{ccccccc}
          &  &  &   1 &  &    &  \\
           &  & 1 &  & 1 &  &    \\ 
          & 1 &  & 0 &  & 1 &    \\
     \end{array} 
     \hspace{1cm}     
     \begin{array}{ccccccccc}
     &     &  &  &   0 &  &    &  \\
     &      &  & 0 &  & 1 &  &    \\ 
      &    & 1 &  & 0 &  & 0 &    \\
     &     1 &  & 0 &  & 1 &  & 0   \\
    0 & & 0 &  & 0 &  & 0 &  & 0   \\
     \end{array} \]
\end{example}
Striker defined and studied these objects while finding a partial bijection between permutation matrices and a subset of TSSCPPs and established the following.
\begin{proposition}
[\protect{\cite[Proposition 2.13]{striker2018permutation}}]
\label{prop:magog}
    The set of all TSSCPP boolean triangles of order $n$ is in bijection with the set of all TSSCPPs in a $2n \times 2n \times 2n$ box.
\end{proposition}

\section{TSSCPP matrices - definition and enumerative properties}
\label{sec:matrices}
In this section, we define a subset of square sign matrices we show is in bijection with the set of TSSCPPs in a $2n\times 2n\times 2n$ box. We then study enumerative properties of these matrices. In the next section, we study the polytope defined as their convex hull.

\subsection{Definition of magog matrices}
\label{subsec:magogdef}
We begin by defining the set of matrices we show in Theorem~\ref{prop:magog_bij} to be in bijection with TSSCPPs.

\begin{definition}\label{ineq:m}
    An $n\times n$ \emph{magog matrix} $A=(a_{ij})$ is a square sign matrix (Definition~\ref{def:sq_sign}) such that:
    \begin{align}          
     \sum_{j'=1}^{j}a_{i+1,j'}+ \sum_{i'=1}^{i+1}a_{i',j+1}-\sum_{i'=1}^{i}a_{i'j} \geq 0 & \text{ for all }
     1 \leq i \leq n-2, 
     1 \leq j \leq n-2.   \label{eq:5magog}
     \end{align}
    
     We refer to \eqref{eq:5magog} as the $(i,j)$-special inequality. Let $\magog_n$ denote the set of $n\times n$ magog matrices. 
\end{definition}

Since magog matrices are square sign matrices, they are characterized as $\{0,1,-1\}$-matrices satisfying \eqref{eq:1magog}-\eqref{eq:4magog} with the added condition of satisfying the $(i,j)$-special inequalities \eqref{eq:5magog}. Also note that the only time an $(i,j)$-special inequality is violated is when the partial sum of column $j$ up to row $i$ is one, but both the partial sum of row $i+1$ up to column $j$ and the partial sum of column $j+1$ up to row $i+1$ are zero.

\begin{example} The matrix below, left, is a non-example of a magog matrix; it violates the $(2,2)$-special inequality because the sum of the green elements is less than the sum of the red elements. The matrix below, right, is not an ASM but all special inequalities are satisfied, so it is a magog matrix. Notice that the $(2,1)$-special inequality is satisfied since the sum of green elements is greater than or equal to the sum of red elements.  

\[\begin{pmatrix} 0 & \color{red}1 & \color{green}0 &\color{black} 0\\ 1 & \color{red}0 & \color{green}0 & \color{black} 0\\ \color{green}0 & \color{green}0 & \color{green} 0 & \color{black} 1\\ 0 & 0 & 1 & 0\end{pmatrix}
\hspace{1in}
\begin{pmatrix}\color{red} 0 & \color{green}0 & 0 & 1\\ \color{red} 0 & \color{green} 1 & 1 & -1\\ \color{green} 1 & \color{green} 0 & 0 & 0\\ 0 & 0 & 0 & 1\end{pmatrix}\]
\end{example}

See Figure~\ref{fig:8sq} for the seven matrices in $\magog_3$ and Figure~\ref{fig:bigFD}, right, for a larger example.

\smallskip
Our first theorem justifies the naming of these matrices by giving a bijection to magog triangles. See Figure~\ref{fig:bigFD} for an example and Remark~\ref{remark:Phi} for some historical notes on the bijection map.

\begin{theorem}
\label{prop:magog_bij}
    The set $\magog_n$ of $n\times n$ magog matrices is in explicit bijection with the set of magog triangles of order $n$, and therefore with the set of all TSSCPPs in a $2n\times 2n\times 2n$ box.
\end{theorem}

\begin{proof}
Given a magog matrix $M\in\magog_n$, we wish to construct a magog triangle $\alpha:=\Psi(M)$.

We first construct an $n\times n$ matrix $M'$ where the entry in the $i$th row, $j$th column is the sum of the entries from rows one through $i$ of the $j$th column of $M$. That is, $M'$ is the column partial sum matrix of $M$. Because the column partial sums of $M$ are nonnegative, as implied by inequality \eqref{eq:3magog}, each of these entries will be either zero or one. Since in addition, the sum of each row or column equals one \eqref{eq:1magog}-\eqref{eq:2magog}, we have that row $i$ of $M'$ has exactly $i$ ones in it.  Because of the $(i,j)$-special inequality~\eqref{eq:5magog}, if the $j$th entry in row $i$ of $M'$ is one, then we will have that there is at least one more one in the first $j+1$ entries of row $i+1$ than in the first $j$ entries of row $i$. 

For the second part of this bijection we take $M'$ as described above and construct a magog triangle $\alpha$ from it. To do this we let row $i$ of $\alpha$ record the indices in row $i$ of $M'$ that contain a one. Thus the rows of $\alpha$ are strictly increasing by construction, yielding~\eqref{eq:tri_rowinc}. Since the $i$th row of $M'$ has exactly $i$ entries equal to one, $\alpha$ will be a triangle, and since each column of $M$ sums to one~\eqref{eq:2magog}, the last row will be $1~2\ldots n$, yielding~\eqref{eq:tribottomrow}. The requirement that the row partial sums of $M$ are nonnegative~\eqref{eq:4magog} together with the $(i,j)$-special inequality~\eqref{eq:5magog} implies via our construction that if $j$ is the $k$th element of row $i$ of $\alpha$, then the $(k+1)$st element of row $i+1$ of $\alpha$ must be less than or equal to $j+1$. This implies \eqref{eq:tri_diagplus1} holds, hence, $\alpha$ is a magog triangle. 

Both steps of the map $\Psi$ are clearly invertible, thus this is a bijection.
By Theorem~\ref{thm:magogtri}, this also establishes a bijection between the set of TSSCPPs in a $2n\times 2n\times 2n$ box and the set $\magog_n$ of magog matrices of order $n$.
\end{proof}

\begin{figure}[htbp]
\begin{center}
\[\begin{array}{ccccl}
      \begin{pmatrix} 0 & 0 & 0 & 1 \\ 0 & 1 & 1 & -1 \\ 1 & 0 & 0 & 0 \\ 0 & 0 & 0 & 1 \end{pmatrix}
     &\longleftrightarrow 
     &\begin{pmatrix}
         0 & 0 & 0 & 1 \\ 0 & 1 & 1 & 0 \\ 1 & 1 & 1 & 0 \\
         1 & 1 & 1 & 1
     \end{pmatrix}
     &\longleftrightarrow 
     &\begin{array}{ccccccc}
          &  &  &   4 &  &    &  \\
           &  & 2 &  & 3 &  &    \\ 
          & 1 &  & 2 &  & 3 &    \\
          1 &  & 2 &  & 3 &  & 4   \\
     \end{array}  
\end{array}\]
\end{center}
\caption{An example of the bijection of Theorem~\ref{prop:magog_bij}. On the left is a magog matrix $M$. In the center is its matrix of partial column sums $M'$. On the right is the corresponding magog triangle $\alpha = \Psi(M)$.}
\label{ex:magog_example}
\end{figure}

 \begin{remark}
 \label{remark:Phi}
     Note that the map to magog triangles used to prove Theorem~\ref{prop:magog_bij} is the same map of~\cite{MILLS1983340} that sends an alternating sign matrix to a \emph{monotone triangle}; an object that is frequently useful in proving theorems about ASMs. 
  
     In fact, the map $\Psi$ is well-defined on any square sign matrix. So given a square sign matrix, one may apply the map $\Psi$ and see whether the resulting triangle is a magog triangle, a monotone triangle, both, or neither to determine whether the matrix is a magog matrix, an alternating sign matrix, both, or neither. We will use both magog matrices and magog triangles throughout the paper, as certain theorems are more directly seen using one object or the other.
     
     The intersection of magog and monotone triangles was studied in \cite{Striker.tetrahedralposet} and also \cite{Ayyer312}, where they were named \emph{gapless} monotone triangles. There is no known counting formula for this intersection; data computed in these papers is given  as sequence A180349 in the OEIS~\cite{oeis}.
 \end{remark}   

In the remaining subsections of Section~\ref{sec:matrices}, we enumerate various subsets of magog matrices corresponding to certain statistic values and give some further conjectures based on enumerative data from SageMath~\cite{sage} computations. 
\subsection{Magog matrices with extreme numbers of negative ones}
\label{subsec:permmagog}
In this subsection, we enumerate magog matrices with no negative ones and then discuss the case of the maximum number of negative ones.

Our first result is that the magog matrices with no negative ones are counted by the Catalan numbers; see the first entry in each row of Table~\ref{tab:MM-1_stats}. This differs from the ASM case, in which all $n!$ permutation matrices are ASMs.

We use the following standard definition of pattern avoidance in permutations.

\begin{definition}
    A permutation $\pi=\pi_{1}\pi_{2}\ldots\pi_{n}\in S_n$ has a $132$ pattern if there are $i,j,k$ such that $i<j<k$ but $\pi_i<\pi_k<\pi_j$. A permutation is \emph{$132$-avoiding} if it has no $132$ pattern.
\end{definition}

Given a permutation $\pi\in S_n$, our convention is that its permutation matrix $(M_{ij})_{i,j=1}^n$ is defined as $M_{ij}=1$ whenever $M_i=j$ and zero otherwise.

\begin{example} The permutation $2431$ has a $132$ pattern. Its permutation matrix is shown below. 
 \hspace{1cm}
        \[\begin{pmatrix}
         0 &1 & 0 & 0 \\ 0&0 & 0 & 1 \\ 0& 0& 1 & 0 \\
         1 & 0 & 0 &0
     \end{pmatrix}\]     
\end{example}

\begin{theorem}
\label{thm:132}
    The magog matrices of order $n$ with no negative ones are the $132$-avoiding permutation matrices.
\end{theorem}

\begin{proof}
    If a magog matrix has no negative ones, then it is a permutation matrix. We will show that the permutation matrix of $\pi$ is a magog matrix if and only if $\pi$ avoids the pattern $132$.  
    
    If $\pi$ is not a $132$-avoiding permutation, then let $a,b,c$ be the values in $\pi$ that form a $132$ pattern. Suppose also that $a$ is the largest value less than $c$ appearing before $b$ in the permutation $\pi$.  We will show the permutation matrix $M$ of $\pi$ is not a magog matrix. 
    
    Suppose the position of $b$ in $\pi$ is $i+1$, so $i+1$ is the row of $M$ that has a one in column $b$. We show the $(i,a)$-special inequality $\eqref{eq:5magog}$ is violated. By construction, $\sum_{i'=1}^{i}M_{i'a}=1$, since $a$ comes before $b$ in $\pi$. Also, $\sum_{j'=1}^{a}M_{i+1,j'}=0$ since $a < b$ and the one in row $i+1$ does not occur until column $b$.   Finally, $\sum_{i'=1}^{i+1}M_{i',a+1}=0$ since $a$ is the greatest column of $M$ less than or equal to $c$ with a one before row $i+1$. So the left-hand side of \eqref{eq:5magog} is $-1$ which is not greater than or equal to zero, so the $(i,a)$-special inequality is violated, showing $M$ is not a magog matrix.

    Now we show that if a permutation matrix $M$ is not a magog matrix, then its corresponding permutation $\pi$ contains a $132$ pattern. Since all permutation matrices satisfy \eqref{eq:1magog}-\eqref{eq:4magog}, we assume $M$ violates the $(i,j)$-special inequality \eqref{eq:5magog}. Since row partial sums are all nonnegative and column partial sums must be zero or one, the only way for $\eqref{eq:5magog}$ to be violated is if $\sum_{i'=1}^{i}M_{i'j}=1$ and $\sum_{j'=1}^{j}M_{i+1,j'}=\sum_{i'=1}^{i+1}M_{i',j+1}=0$. The equality $\sum_{j'=1}^{j}M_{i+1,j'}=0$ implies $\sum_{j'=j+2}^{n}M_{i+1,j'}=1$; letting $b$ be the column of the one in row $i+1$, we see $b\geq j+1$. Finally, $\sum_{i'=1}^{i+1}M_{i',j+1}=0$ implies $\sum_{i'=i+2}^{n}M_{i',j+1}=1$, so the one in column $j+1$ must be in row $i+2$ or later. Thus, 
    $j < j+1 < b$ and these values form a $132$ pattern in $\pi$. 
\end{proof}
\begin{example}Using the matrix for $2431$ we see that the $(1,2)$-special inequality is violated, as the sum of the green elements is less than the red element. 
     \hspace{1cm}
        \[\begin{pmatrix}
         0 &\color{red}1 &\color{green} 0 & 0 \\ \color{green}0&\color{green}0 & \color{green}0 & 1 \\ 0& 0& 1 & 0 \\
         1 & 0 & 0 &0
     \end{pmatrix}\]     
\end{example}
It is well-known (see e.g.\ \cite[Section 2.2.1, Exercise 4]{knuth97}) that the set of $132$-avoiding permutations of $n$ is enumerated by the $n$th Catalan number.
\begin{corollary}
    The number of $n\times n$ magog matrices with no negative ones equals the $n$th Catalan number $\displaystyle\frac{1}{n+1}\displaystyle\binom{2n}{n}$.
\end{corollary}

We now turn to the other extreme and study magog matrices with the maximum number of negative ones. We show this maximum number is the same as for alternating sign matrices.
Let $\max_{-1}(X)$ be the maximum number of negative ones in a matrix over a set of matrices $X$.

\begin{lemma}
\label{thm:max-1}
    The maximum number of negative ones in any $n\times n$ square sign matrix is given by $\max_{-1}(\sign_n)=\lfloor\frac{n-1}{2}\rfloor\lceil\frac{n-1}{2}\rceil$.
\end{lemma}
\begin{proof}
    First we show that  $\max_{-1}(\sign_n) \ge \lfloor\frac{n-1}{2}\rfloor\lceil\frac{n-1}{2}\rceil$ by constructing an $n\times n$ sign matrix $M_n$ with $\lfloor\frac{n-1}{2}\rfloor\lceil\frac{n-1}{2}\rceil$ negative ones. 
    Let row $\lfloor\frac{n+1}{2}\rfloor$ of $M_n$ alternate between $1$ and $-1$ as in: $1-1~~1\cdots 1-1~~1$ if $n$ is odd, and append a $0$ if $n$ is even. The row directly above this row is given as $0~~1-1\cdots -1~~1~~0$ if $n$ is odd, and with a zero appended if $n$ is even. Each successive row above this is likewise alternating with one more $0$ at the beginning and end, so that row $i$ for $1\leq i\leq \lfloor\frac{n+1}{2}\rfloor$ has $\lfloor\frac{n+1}{2}\rfloor-i$ zeros at the beginning and end with alternating ones and negative ones in the middle when $n$ is odd, and a zero appended if $n$ is even. So there are $i-1$ negative ones in row $i$ for $1\leq i\leq\lfloor\frac{n+1}{2}\rfloor$.     
    The bottom half of the matrix is constructed to make a matrix that is symmetric via a half-turn rotation. See Example~\ref{ex:max-1ex}.      One may easily show by direct computation that the matrix $M_n$ has the required number of negative ones. (Note that when $n$ is even, there is another matrix with this many negative ones constructed as the vertical reflection of $M_n$.)
    
    Now we show that $\max_{-1}(\sign_n) \le \lfloor\frac{n-1}{2}\rfloor\lceil\frac{n-1}{2}\rceil$. 
    Suppose by way of contradiction that some $B=(b_{ij}) \in \sign_n$ has more negative ones than the matrix $A$ described above. Then there must be some row $i$ where $B$ has more negative ones than $A$; suppose for simplicity this row is in the top half of the matrix, so $1\leq i\leq \lfloor\frac{n+1}{2}\rfloor$.  Now since $A$ has $i-1$ negative ones in row $i$, then $B$ has at least $i$. Now let $J$ be the set of columns in $B$ with a negative one in row $i$; the partial sum of each of these columns from the top to row $i-1$ must equal one. So $i = \sum_{j \in J}\sum_{i'=1}^{i-1}b_{i'j} \le \sum_{j' =1}^{n}\sum_{i'=1}^{i-1}b_{i'j'}=i-1$ since each full row sums to one. This is a contradiction. 
\end{proof}
\begin{example}
\label{ex:max-1ex}
$M_5$ is an example of the square sign matrix with the maximum number of negative ones when $n$ is odd. $M_6$ is an example when $n$ is even. 
    \[\begin{array}{ccl}
      M_5=\begin{pmatrix} 0 & 0 & 1 & 0 & 0 \\ 0 & 1 & -1 & 1 & 0\\ 1 & -1 & 1 & -1 & 1 \\ 0 & 1 & -1 & 1 & 0 \\ 0 & 0 & 1 & 0& 0 \end{pmatrix}

    \hspace{1cm}
       M_6=\begin{pmatrix} 0 & 0 & 1 & 0 & 0 & 0 \\ 0 & 1 & -1 & 1 & 0 & 0\\ 1 & -1 & 1 & -1 & 1 & 0\\0 & 1 & -1 & 1 & -1 & 1\\ 0 & 0 & 1 & -1 & 1 & 0 \\ 0 & 0 & 0 & 1 & 0& 0 \end{pmatrix}
\end{array}\]
\end{example}
\begin{corollary}
    The maximum number of negative ones in any $n\times n$ ASM or magog matrix is given by $\max_{-1}(ASM_{n})=\max_{-1}(\magog_{n}) =\lfloor\frac{n-1}{2}\rfloor\lceil\frac{n-1}{2}\rceil$.
\end{corollary}

\begin{proof}
    It is easy to see that the matrix $M_n=(m_{ij})$ constructed in the proof of Lemma~\ref{thm:max-1} is an ASM, as the nonzero elements along any rows or column clearly alternate between positive and negative entries. 

    To show it is also a magog matrix, we must show that the $(i,j)$-special inequality holds. It is helpful to look at $M_5$ and $M_6$ of Example~\ref{ex:max-1ex}. If $m_{ij}=-1$ there is nothing to prove. If $m_{ij}=0$ and $i \leq \lfloor\frac{n+1}{2}\rfloor$, the column $j$ partial sum will be zero, so we are done.

    If $m_{ij}=0$ and $i > \lfloor\frac{n+1}{2}\rfloor$, the column $j$ partial sum is one but the row $i + 1$ partial sum is one and the column $j + 1$ partial sum is either zero or one. If $m_{ij}=1$, then $m_{i+1,j+1}$ is zero or one. If it is one, then the inequality is satisfied. If it is zero, that means both row $i+1$ and column $j+1$ have all their nonzero elements above and to the left of $m_{ij}$, and the inequality is satisfied.
\end{proof}

\begin{remark}
    Note the number of $n\times n$ ASMs with the maximum number of negative ones is one if $n$ is odd and two if $n$ is even, as seen in Table~\ref{tab:asmnegone}, since the matrices in the above construction are the only such ASMs. See \cite[Sec.~2.1]{behrend.multiply.refined} for further discussion. However, the number of $n\times n$ magog matrices with the maximum number of negative ones is much larger, as seen in Table~\ref{tab:MM-1_stats}. 
\end{remark}
\begin{table}[htbp]
\centering
\begin{tabular}{|c | c |}
 \hline
 $n$  & Distribution of negative ones in magog matrices \\
 \hline\hline
 3   & 5,2  \\ 
 \hline
 4 & 14,21,7  \\
 \hline
 5  & 42,149,166,64,8 \\
 \hline
 6  & 132, 892, 2186, 2424, 1373, 379,50 \\
 \hline
 7& 429, 4857, 21567, 48323, 62565, 48933, 23684, 6836, 1075,
79\\
 \hline
\end{tabular}
\caption{Data on negative ones in magog matrices.}
    \label{tab:MM-1_stats}
\begin{tabular}{|c | c | c |} 
 \hline
 $n$  & Distribution of negative ones in ASMs \\
 \hline\hline
 3  & 6,1   \\ 
 \hline
 4 & 24,16,2   \\
 \hline
 5  & 120,200,94,14,1  \\
 \hline
 6  & 720, 2400, 2684, 1284, 310, 36, 2 \\
 \hline
 7 & 5040, 29400, 63308, 66158, 38390, 13037, 2660, 328, 26, 1\\
\hline
\end{tabular}
\caption{For comparison: data on negative ones in ASMs.}
    \label{tab:asmnegone}

\end{table}

\subsection{Magog matrices with boundary ones in specified entries} 
\label{subsec:bdrymagog}
For alternating sign matrices, the existence of a one in any of the four corners implies the rest of the entries in its row and column must be zero. Furthermore, removing the entire row and column results in an ASM one size smaller.  In this subsection, we will see this holds for the bottom left corner of a magog matrix, but not for any of the other corners. Moreover, while it is true for magog matrices (and in fact all sign matrices) that the top row, left column, and bottom row must consist of all zeros except a single one, this is not true of the last column, as seen in Figure~\ref{fig:8sq}. 

Let $\magog_{n}(i,j)$ be the set of magog matrices of order $n$ with a one in row $i$ column $j$. See Table~\ref{tab:MMfirstone_stats} for data on the cardinality of these sets for $(i,j)$ in the first row, first column, or last row; see also Table~\ref{tab:ASMfirstone_stats} for the analogous ASM data, for comparison.

The following is the main theorem of this subsection.
\begin{theorem} 
\label{thm:bdry_ones}
The following enumerative identities hold:
    \begin{enumerate}
        \item \label{310_1} For all $n\geq 1$,  $|\magog_{n}(1,1)|=1$,
        \item  \label{310_2}   For all $n>1$, $|\magog_{n}(n,1)|=|\magog_{n}(n,2)|=|\magog_{n-1}|=\displaystyle\prod_{j=0}^{n-2}\frac{(3n-2)!}{(n-1+j)!}$,
        \item\label{310_3} For all $n>1$, $|\magog_{n}(1,n)|=|\magog_{n}(1,n-1)|$,
        \item\label{310_4} For $n\geq 1$, $|\magog_{n}(2,1)|=\displaystyle\frac{1}{n+1}\binom{2n}{n}-1$, and
        \item\label{310_5} For all $n\geq 1$,
    $|\magog_{n}(1,2)|=2^{n-1}-1$.
    \end{enumerate}
\end{theorem}

We prove this theorem via the following lemmas. The first shows there is only one magog matrix with a one in the upper left corner.
\begin{lemma}
\label{upperleft1}
The only matrix in $\magog_n$ with a $1$ in the upper left corner is the identity matrix.
\end{lemma}
\begin{proof}
    Let $M\in \magog_{n}(1,1)$, so its upper left corner entry equals $1$.  Let $\alpha = (\alpha_{ij})$ be its corresponding magog triangle of Definition~\ref{def:magogtri}, given by performing the inverse map $\Psi^{-1}$ from the proof of Theorem~\ref{prop:magog_bij}. Then the top entry of $\alpha$ equals $1$, and thus by~\eqref{eq:tri_diagplus1}, the last northwest to southeast diagonal is $1~2~\ldots~n$. The defining magog triangle inequalities \eqref{eq:tri_rowinc} and \eqref{eq:tri_diagplus1}  imply \eqref{eq:tri_implied}, which says that the southwest to northeast diagonals are weakly increasing. Since the last row of the triangle must be $1~2~3\ldots n$ by \eqref{eq:tribottomrow}, this shows that all the entries of the $i$th southwest to northeast diagonal of $\alpha$ must equal $i$. The magog matrix corresponding to this magog triangle via the map $\Psi$ is the identity matrix. 
\end{proof}

The next lemma on the enumeration of these sets is the following, which shows the number of $n\times n$ magog matrices with a one in the last row and first column equals the total number of $(n-1)\times (n-1)$ magog matrices, whose enumeration is given by the product formula \eqref{asmform}. This refined enumeration also holds for alternating sign matrices. The subsequent lemma shows that unlike in the case of ASMs, this is also the number of magog matrices whose last row one is in column two. 

\begin{lemma}
\label{thm:lastrowfirstcol}
    For all $n>1$, $|\magog_{n}(n,1)|=|\magog_{n-1}|.$
\end{lemma}
\begin{proof}
    We give a bijection from $\magog_{n}(n,1)$ to the set of all $(n-1)\times (n-1)$ magog matrices $\magog_{n-1}$. Note that when the last row has a one in column one, the rest of the entries in the first column and the rest of the entries in the last row are all zero. We may remove the first column and last row to get an $(n-1)\times (n-1)$ magog matrix. Likewise, adding a row at the bottom and a column to the left with a one at their intersection to any $(n-1)\times (n-1)$ magog matrix does not violate any $(i,j)$-special inequalities \eqref{eq:5magog}, and so will give us a new $n\times n$ magog matrix. 
    \end{proof}

 \begin{lemma}
 \label{n1n2}
    For all $n>1$, $|\magog_{n}(n,1)|=|\magog_{n}(n,2)|$.
 \end{lemma}
 \begin{proof}
    For $n>1$, we give a bijection between 
    $\magog_{n}(n,1)$ and $\magog_{n}(n,2)$. 
       Choose $A=(a_{ij}) \in \magog_{n}(n,1)$, meaning $a_{n1}=1$.  Let $\alpha = (\alpha_{ij})$ be its corresponding magog triangle given by performing the inverse map $\Psi^{-1}$ from the proof of Theorem~\ref{prop:magog_bij}. If the unique $1$ in the last row of $A$ is in column $j$, this means $j$ is the unique number missing from row $n-1$ of $\alpha$.
    Since $a_{n1}=1$, row $n-1$ of $\alpha$ is exactly $2~3~4\ldots n$. We may change this first entry $\alpha_{n-1,1}=2$ to a $1$ without violating any of the defining inequalities \eqref{eq:tri_len}-\eqref{eq:tribottomrow}, as the only constraints are that $\alpha_{n-1,1}+1\geq \alpha_{n,1}=2$, and that $\alpha_{n-1,1}<\alpha_{n-1,2}=3$. Thus, this changes row $n-1$ to $1~3~4\ldots n$, which maps by $\Psi$ to a matrix in $\magog_{n}(n,2)$. Conversely, if $A$ were in  $\magog_{n}(n,2)$, its magog triangle would have row $n-1$ equal to $1~3~4\ldots n$, and we could similarly change the $1$ to a $2$, producing a magog triangle whose magog matrix is in $\magog_{n}(n,1)$. So $|\magog_{n}(n,1)|=|\magog_{n}(n,2)|$.
\end{proof}

We now show that the number of magog matrices whose top row one is in column $n$ equals the number whose top row one is in column $n-1$. Note from Table~\ref{tab:MMfirstone_stats} that the cardinality of these sets grows quickly; we do not have even a conjectured enumeration formula.

\begin{lemma}
\label{thm:firstrowones}
    For all $n>1$, $|\magog_{n}(1,n)|=|\magog_{n}(1,n-1)|$. 
\end{lemma}
\begin{proof}
Given a magog triangle, the only constraints on the top entry are that it is at most $n$ and greater than or equal to its southeast diagonal neighbor minus one. If the top entry of a magog triangle is $n$ or $n-1$; this is always greater than or equal to its southeast diagonal neighbor minus one, since by \eqref{eq:tri_len}, all entries in a magog triangle are at most $n$. Given a magog matrix in  $\magog_{n}(1,n)$, its magog triangle has top entry $n$. We may change this to $n-1$ without violating the defining inequalities. Likewise, given a magog matrix in  $\magog_{n}(1,n-1)$, its magog triangle has top entry $n-1$. We may change this to $n$ without violating the defining inequalities. Thus we have a bijection between $\magog_{n}(1,n)$ and $\magog_{n}(1,n-1)$ given by changing the top entry of the corresponding magog triangle from $n$ to $n-1$. Therefore, $|\magog_{n}(1,n)|=|\magog_{n}(1,n-1)|$.
\end{proof}

The next lemma enumerates magog matrices whose first column $1$ is in the first or second row.
\begin{lemma}
\label{thm:Catalan}
    The number of $n\times n$ magog matrices with a one in the first or second entry of its first column is given by the $n$-th Catalan number $\displaystyle\frac{1}{n+1}\binom{2n}{n}$.
\end{lemma}
\begin{proof}

    Given $A\in\magog_n$ with a one in the first or second entry of its first column, its corresponding magog triangle $\alpha$ has a $1$ as the first entry of the second row. If we consider the triangle formed by this entry and all entries below it, this must be a magog triangle of $n-1$ rows. By Lemma~\ref{upperleft1}, this must be the magog triangle whose corresponding magog matrix is the $(n-1)\times(n-1)$ identity matrix. Thus, all entries of $\alpha$ except the last north to southeast diagonal are fixed to be the lowest number possible, and thus place no constraints on that last diagonal. 

    Consider the last diagonal of $\alpha$: $\alpha_{1,n-1}~ \alpha_{2,n-1}~ \ldots~\alpha_{n-1,n-1}~n$. By the magog triangle inequality \eqref{eq:tri_diagplus1}, we have $\alpha_{ij}+1\geq \alpha_{i+1,j}$.  Thus we have $\alpha_{n-1,n-1}+1\geq n$ so $\alpha_{n-1,n-1}\geq n-1$. By iterating this reasoning, $\alpha_{i,n-1}\geq i$.  
   Then
    if we let $y_i=\alpha_{i,n-1}-i$, we obtain a sequence $y_1, y_2,\ldots, y_n$ with $y_n=0$, $0\leq y_i\leq n-i$, and $y_i\geq y_{i+1}$. Such sequences are known to be counted by the  Catalan numbers $\frac{1}{n+1}\binom{2n}{n}$ (see e.g.\cite[Ex.~6.19(s)]{EC2}). 
\end{proof}

The final lemma enumerates magog matrices whose top row $1$ is in column one or two.
\begin{lemma}
\label{thm:1221}
The number of $n\times n$ magog matrices with a one in the first or second entry of its first row is given by $2^{n-1}$.
\end{lemma}
\begin{proof}
   Given a magog matrix $A\in\magog_n$ whose top row $1$ is in column one or two, the top entry of the corresponding magog triangle $\alpha$ is $1$ or $2$. We wish to show the number of such magog triangles if $2^{n-1}$.
   We prove this by induction. The case $n=1$ is clear. Now suppose there are $2^{n-2}$ magog triangles of order $n-1$ and top entry $1$ or $2$; call this set $\mathcal{A}$. For each $\alpha\in\mathcal{A}$, we construct two magog triangles of order $n$. The first is constructed by appending a new final row $1~2~\ldots~n$ to $\alpha$; this row is clearly strictly increasing, so it satisfies~\eqref{eq:tri_rowinc}. Since the last row of $\alpha$ is $1~2~\ldots n-1$, the inequality~\eqref{eq:tri_diagplus1} is also clearly satisfied. The second new magog triangle is constructed from $\alpha$ by adding to $\alpha$ a new north to southeast diagonal $2~3~4~\ldots~n-1~n~n$. This diagonal itself satisfies \eqref{eq:tri_diagplus1}. Also, the last diagonal of $\alpha$ is at most $2~3~\ldots n-2~n-1~n-1$ by \eqref{eq:tri_diagplus1}, so \eqref{eq:tri_rowinc} is satisfied. Thus, there are at least $2(2^{n-2})=2^{n-1}$ magog triangles of order $n$ with top entry $1$ or $2$. 
   Suppose there was a magog triangle $\beta$ of order $n$ with top entry at most $2$ that did not result from this construction. The last entry in row $n-1$ is either $n-1$ or $n$. If it is $n$, then the last diagonal is $2~3~\ldots~n-1~n~n$, otherwise we could not have top entry at most $2$. The number $n$ does not appear in any entry outside this diagonal, so the rest of the triangle without this diagonal is a magog triangle of order $n-1$, implying $\beta$ did come from this construction. If the last entry in row $n$ of $\beta$ is $n-1$, then row $n-1$ is $1~2~\ldots n-1$ and the rest of the triangle without the bottom row is a magog triangle of order $n-1$ with top entry $1$ or $2$, so $\beta$ came from this construction in this case also. So, there are exactly $2^{n-1}$ magog triangles of order $n$ with top entry at most $2$. 
\end{proof}

\begin{proof}[Proof of Theorem~\ref{thm:bdry_ones}]
Lemma~\ref{upperleft1} shows $|\magog_n(1,1)=1|$, and thus proves Part (\ref{310_1}). Lemmas~\ref{thm:firstrowones} and \ref{n1n2} prove the first two equalities of (\ref{310_2}). The last equality follows by Theorem~\ref{prop:magog_bij} and Andrews' proof that \eqref{asmform} counts TSSCPPs in a $2n\times 2n\times 2n$ box~\cite{ANDREWS199428}. Lemma~\ref{thm:lastrowfirstcol} proves (\ref{310_3}). Lemma~\ref{thm:Catalan} shows $|\magog_n(1,1)|+|\magog_n(2,1)|$ equals the $n$th Catalan number, so combined with Lemma~\ref{upperleft1}, we have (\ref{310_4}). Finally, Lemma~\ref{thm:1221} shows $|\magog_n(1,1)|+|\magog_n(1,2)|$ equals $2^n$, so combined with Lemma~\ref{upperleft1}, we have (\ref{310_5}). 
\end{proof}

We end with a remark comparing Theorem~\ref{thm:bdry_ones} with analogous enumerations for alternating sign matrices; see Tables~\ref{tab:MMfirstone_stats} and \ref{tab:ASMfirstone_stats}.
\begin{remark}
    The only part of Theorem~\ref{thm:bdry_ones} whose analogous statement also holds for ASMs is the part of (2) from Lemma~\ref{thm:lastrowfirstcol}: that the number of matrices with a one in the lower left corner equals the number of matrices of one size smaller. There is a counting formula for the number of ASMs whose $1$ in the top row is in column $j$; this formula is a nice product formula, conjectured in \cite{MILLS1983340} and proved in \cite{Zeilberger_refined1996}.
\end{remark}

\begin{table}[htbp]
\centering
\begin{tabular}{|c | c | c | c |} 
\hline 
 $n$   & Distribution of & Distribution of  & Distribution of \\
   &   first row one & first column one & last row one \\
 \hline\hline
 3   & 1,3,3 & 1,4,2 & 2,2,3 \\ 
 \hline
 4  & 1,7,17,17 & 1, 13, 21, 7 & 7,7,12,16  \\
 \hline
 5  & 1,15,75,169,169&1, 41, 177, 168, 42&  42, 42, 77, 119, 149 \\
 \hline
 6   & 1, 31, 304, 1328, 2886, 2886 & 1, 131, 1462, 3268, 2145, 429& 429, 429, 816, 1380, 1988, 2394\\
 \hline
 7&   1, 63, 1190, 9690, 39444, & 1, 428, 12506, 63570,& 7436, 7436, 14443, 25883,   \\ 
   & 83980, 83980 & 89791, 44616, 7436 & 41028, 56212, 65910\\
   \hline
\end{tabular}
\caption{Distribution of first row, last row, and first column ones for magog matrices.}
    \label{tab:MMfirstone_stats}
\end{table}    

\begin{table}[htbp]
\begin{tabular}{|c | c | c |} 
 \hline
 $n$   & Distribution of \\
   &  first row / last row / first column one \\
 \hline\hline
 3   & 2,3,2  \\ 
 \hline
 4  & 7,14,14,7  \\
 \hline
 5  & 42,105,135,105,42  \\
 \hline
 6   & 429, 1287, 2002, 2002, 1287, 429 \\
 \hline
 7 & 7436, 26026, 47320, 56784, 47320, 26026, 7436\\
\hline
\end{tabular}
\caption{Distribution of first row/ last row/ and first column ones of alternating sign matrices for comparison.}
    \label{tab:ASMfirstone_stats}
\end{table}

\subsection{Magog matrices with a specified number of inversions}
\label{subsec:magoginv}

We now give some enumerations and conjectures related to inversions and positive inversions

\begin{definition}
    For any square sign matrix $A=(a_{ij})$, we define the \emph{inversion number} $\inv(A)$ to be
        \[
\inv(A)=\sum_{1\leq k < i \leq n}\sum_{1\leq j<l\leq n}a_{ij}a_{kl}.\]
    The positive inversion number $\posinv(A)$ of a square sign matrix is
    \[
        \posinv(A)=\inv(A)-\mathcal{N}(A)
    \]
    where $\mathcal{N}(A)$ equals the number of negative one entries in $A$. We also use the following notation for fixed $1\leq k,l \leq n$:
\[
        \inv_{k,l}(A)=\sum_{k<i\leq n}\sum_{1\leq j<l}a_{kl}a_{ij}.
    \]
\end{definition}
The refined statistic $\inv_{k,l}(A)$ can be thought of as the number of inversions involving $a_{kl}$ and entries southwest of it. Note that $\inv(A)=\sum_{l=1}^{n}\sum_{k=1}^{n} \inv_{k,l}(A)$.

For any permutation  $\pi$, $\posinv(\pi)=\inv(\pi)$, which equals the usual inversion number of $\pi$ as a permutation. 
Inversions and positive inversions of alternating sign matrices were first studied in \cite{MILLS1983340} and \cite{ROBBINS1986169}. Statistics for inversions and positive inversions of ASMs and magog matrices are contained in Tables \ref{tab:MMinversion_stats} and \ref{tab:ASMinversion_stats}. For more information on inversions and positive inversions see 
\cite[Sections 2.1 and 5.2]{behrend.multiply.refined}.

We begin with the following lemma.
\begin{lemma}
\label{lem:maxinv}
    The number of inversions 
    of a square sign matrix is at most $\binom{n}{2}$. Moreover, the only matrix with $\binom{n}{2}$ inversions is the antidiagonal permutation matrix, that is, the matrix with ones down the main antidiagonal and zeros elsewhere.
\end{lemma}
\begin{proof}
    Choose a square sign matrix $A=(a_{ij})$.
    Fix $k$ and $l$ between $1$ and $n$. Observe that $\inv_{k,l}(A)=a_{kl}\sum_{k<i\leq n}\sum_{1\leq j<l}a_{ij}=a_{kl}\sum_{1\leq j<l}\sum_{k<i\leq n}a_{ij}\leq a_{kl}\sum_{1\leq j<l}1=a_{kl}(l-1)\leq l-1$, where the first inequality holds because partial column sums from the bottom must be at most $1$ since partial column sums from the top are at least zero and each column sums to $1$.
    Now consider $I=(k_1,\ldots,k_m)$ to be the set of all rows whose entry in column $l$ is nonzero. For $r\geq 1$ we know that since $a_{k_{2r},l}=-1$ we have $\inv_{k_{2r},l}(A)+\inv_{k_{2r+1},l}(A)=(-1)\sum_{i=k_{2r}+1}^{k_{2r+1}}\sum_{j<l}a_{ij}$ which is at most zero since this is a sum of row partial sums. This means $\sum_{k=1}^{n}\inv_{k,l}(A) = \inv_{k_1,l}(A)+\sum_{r=1}^{\frac{m}{2}-1}\inv_{k_{2r},l}(A)+\inv_{k_{2r+1},l}(A) \leq  \inv_{k_1,l}(A) \leq l-1$. Since the total number of inversions is the sum of inversions for each column,  $\inv(A) \leq \sum_{l=1}^{n}\sum_{k=1}^{n} \inv_{k,l}(A) \leq 0+1+2+3+...+n-1 = \binom{n}{2}$.
    
    It is evident that the $n\times n$ antidiagonal permutation matrix has $\binom{n}{2}$ inversions. It remains to show that no other square sign matrix has $\binom{n}{2}$ inversions. It is well-known that any other permutation matrix has fewer inversions, since any transposition applied to the reverse permutation $n~\cdots~3~2~1$ decreases the number of inversions.

  So we need only show that any square sign matrix $M=(m_{ij})$ with a negative one has fewer than $\binom{n}{2}$ inversions. Let $i'$ be the first row in $M$ with a negative one, $j'$  be the column of the first one in row $i'$, and $j''$ be the column of the second one in row $i'$. Since $i'$ is the first row with a negative one, we know the one in column $j''$ is the first nonzero entry in its column. Note $\sum_{i'< i\leq n}m_{i',j''}m_{i,j'}=0$ because $m_{i'j'}=1$, so the partial sum of column $j'$ from row $i'+1$ to row $n$ must be zero. Thus,
    \[\displaystyle\inv_{i',j''}(M)=\sum_{i'< i\leq n}(m_{i',j''}m_{i,j'}+\displaystyle\sum_{\substack{1 \leq j <j''\\j\neq j'}}m_{i',j''}m_{ij})=\sum_{i'< i\leq n}\displaystyle\sum_{\substack{1 \leq j <j''\\j\neq j'}}m_{i',j''}m_{ij}\leq j''-2.\]
By the reasoning in the prior paragraph, $\sum_{i=1}^n\inv_{i,j''}(M)\leq \inv_{i',j''}(M)\leq j''-2$, so \[\inv(M)\leq \sum_{i=1}^{n}\sum_{j=1}^{n} \inv_{i,j}(M) \leq j''-2 +\sum_{\substack{1 \leq j \leq n\\j\neq j''}}(j-1) \leq \binom{n}{2}-1\]
completing the proof.
\end{proof}

\begin{corollary}
    The number of inversions 
    of a magog matrix is at most $\binom{n}{2}$.
\end{corollary}
Note that this bound is tight (as shown in the next theorem), since the antidiagonal permutation matrix, which has $\binom{n}{2}$ inversions, is a magog matrix.

We now count the magog matrices with certain numbers of inversions. Let $\magog_n(k~\inv)$ denote the set of $n\times n$ magog matrices with inversion number $k$. Let $\magog_n(k~\posinv)$ denote the set of $n\times n$ magog matrices with positive inversion number $k$.  See Table~\ref{tab:MMinversion_stats} for data on the cardinality of these sets for $(i,j)$ in the first row, first column, or last row; see also Table~\ref{tab:ASMinversion_stats} for the analogous ASM data, for comparison.

The following is the main theorem of this subsection.
\begin{theorem}
\label{thm:inv}
The following enumerations hold:
\begin{enumerate}
    \item  For all $n\geq 1$, $|\magog_n(0~\inv)|=|\magog_n(0~\posinv)|=1$,
    \item For all $n\geq 2$, $|\magog_n(1~\inv)|=1$, 
    \item For all $n\geq 2$, $|\magog_n(\binom{n}{2}-1~\posinv)|=n-1$,
    \item For all $n\geq 2$, $|\magog_n(\binom{n}{2}~\inv)|=|\magog_n(\binom{n}{2}~\posinv)|=1$, and
    \item For all $n\geq 3$, $|\magog_n(2~\inv)|=n+1$.
\end{enumerate}
\end{theorem}
\begin{proof}
    There is a unique square sign matrix with zero inversions: the matrix of the identity  permutation $1~2~3\ldots n$. So $|\magog_n(0~\inv)|=1$. Each $-1$ will add a minimum of two inversions, so the identity is also the only magog matrix with zero positive inversions. So $|\magog_n(0~\posinv)|=1$. Thus $(1)$ is proved.
    
    While there are $n-1$ permutation matrices with one inversion, constructed by inverting any pair $(i, i+1)$ from $1~2~3\ldots n$, the only one whose matrix is a magog matrix is $2~1~3~4\ldots n$. No matrix for a permutation $1~2\ldots i-1~i+1~i~i+2\ldots n$  is a magog matrix, since the $i-1,i+1,i$ would form a $132$-pattern, which is forbidden, by Theorem~\ref{thm:132}. Any square sign matrix with a negative one would have at least two inversions. Thus, $\magog_n(1~\inv)=1$, proving $(2)$.
    
    By Lemma~\ref{lem:maxinv}, the maximum number of inversions of a square sign matrix is $\binom{n}{2}$ and the unique matrix achieving this is the matrix with ones down the main antidiagonal and zeros elsewhere. This matrix is a magog matrix, since it satisfies \eqref{eq:5magog}. Thus, $|\magog_n(\binom{n}{2}~\inv)|=|\magog_n(\binom{n}{2}~\posinv)|=1$, yielding $(4)$.
    
    For $n\geq 2$, the number of magog matrices with $\binom{n}{2}-1$ positive inversions is $n-1$; these are constructed from this antidiagonal matrix by uninverting any pair of antidiagonally adjacent ones. So $|\magog_n(\binom{n}{2}-1~\posinv)|=n-1$. Note there are no magog matrix with a negative one that has this many positive inversions, since then it would have at least $\binom{n}{2}$ inversions, which is only true of the antidiagonal permutation matrix. Thus, we have shown $(3)$.

        We now show $(5)$, that the number of magog matrices with exactly two inversions is $n+1$ for $n \geq 3$. One may directly compute that the set $X_3$ comprised of the following matrices is the set of $3\times 3$ magog matrices with  inversion number  two: 
    \[\left(\begin{array}{cccc}
         0 & 0 & 1   \\
         1 & 0 & 0 \\
         0 & 1 & 0  \\         
    \end{array}\right)
    \left(\begin{array}{cccc}
         0 & 1 & 0   \\
         0 & 0 & 1 \\
         1 & 0 & 0  \\         
    \end{array}\right)
    \left(\begin{array}{cccc}
         0 & 1 & 0   \\
         1 & -1 & 1 \\
         0 & 1 & 0  \\         
    \end{array}\right)
    \left(\begin{array}{rrr}
         0 & 0 & 1   \\
         1 & 1 & -1 \\
         0 & 0 & 1  \\         
    \end{array}\right).\] 
    
    We follow a recursive procedure to construct $X_n$ from $X_{n-1}$ for arbitrary $n>3$ and prove $X_n$ has cardinality $n+1$ inductively.
    To each matrix in $X_{n-1}$,  we add a row and a column at the end with all zeros except a single one in the lower right corner; this yields $n$ such matrices. We create one more matrix by taking the single matrix in $X_{n-1}$ that has three nonzero entries in the last column and add a row between rows $n-2$ and $n-1$ and a column between columns $n-2$ and $n-1$ with a one at their intersection and zeros at all other entries; an example is shown below for $n=4$. 
    \[\left(\begin{array}{cccc}
     0&0 & 0 & 1   \\
     1&1 & 0 & -1 \\
     0&0 & 1 & 0  \\  
     0& 0&0&1
  \end{array}\right)\]  
    Thus, $X_n$ has $n+1$ elements by induction, all of which are magog matrices and have inversion number two by construction. 
    
    All other sign matrices not resulting from this construction either have more inversions or violate an $(i,j)$-special inequality. 
    We show this by induction. Fix $n>3$ and suppose the $n$ matrices resulting from this construction are all the $(n-1)\times (n-1)$ magog matrices with an inversion number of two. Suppose there is a matrix $A=(a_{ij})$ in $\magog_n$ not resulting from this construction that also only has two inversions. It cannot have a single one at the bottom right entry and all zeros otherwise in the last row and last column, since otherwise the $(n-1)\times (n-1)$ matrix resulting from removing the last row and column would be a magog matrix of size one less with two inversions, so our matrix would have come from the construction. Recall also that each $-1$ will add a minimum of two inversions, so we may have at most one $-1$. 
    
    Now, we  show the last row of $A$ must be all zeros followed by a one. Suppose this is not the case. Then either we must have ones at $a_{n,n-1}$ and $a_{n-1,n}$ or ones at $a_{n,n-2}$, $a_{n-2,n-1}$ and $a_{n-1,n}$, otherwise there will be too many inversions. If there are no negative ones, then the $(n-2,n-2)$-special inequality and $(n-2,n-3)$-special inequality will be violated, respectively. If we do have a negative one, the only magog matrices with a single $-1$  and these specified ones, respectively, are shown below (for $n=5$). Both have many inversions. Thus the last row one must be in column $n$.
    \[\begin{array}{ccl}
      \begin{pmatrix} 0 & 1 & 0 & 0 & 0 \\ 0 & 0 & 1 & 0 & 0\\ 0 & 0 & 0 & 1 & 0 \\ 1 & 0 & 0 & -1 & 1 \\ 0 & 0 & 0 & 1& 0 \end{pmatrix}

    \hspace{1cm}
       \begin{pmatrix}  0 & 1 & 0 & 0 & 0 \\ 0 & 0 & 1 & 0 & 0 \\ 0 & 0 & 0 & 1 & 0\\1 & 0 & -1 & 0 & 1 \\ 0 & 0 & 1 & 0 & 0   \end{pmatrix}
\end{array}\]

    If the last column has no negative ones, then it also may have no other one, and $A$ is a matrix from the construction.
    So we suppose the last column has a single $-1$ in row $k_2$ and a one above it in row $k_1$. Note that since the $-1$ is in the last column there must be two ones preceding it in its row. Since both of those ones will count towards  $\inv_{k_1,n}(A)$ but not $\inv_{k_2,n}(A)$, we have that $\inv_{k_1,n}(A) \geq \inv_{k_2,n}(A)+2$ with equality if and only if $k_1=k_2-1$. Finally, notice that the ones from all other columns must occur in descending order to avoid any extra inversions. From here it follows that  $k_1=1$ because if it was anything else, the $(k_1-1,k_1)$-special inequality would not hold. 
\end{proof}

We end with a remark comparing Theorem~\ref{thm:inv} with the corresponding data on ASMs; see Tables~\ref{tab:MMinversion_stats} and \ref{tab:ASMinversion_stats}.
\begin{remark}
    Note that the same enumerations as Theorem~\ref{thm:inv} (1), (3), and (4) hold for alternating sign matrices (see \cite[Sec.~2.1]{behrend.multiply.refined} for further discussion).  Indeed, all the matrices in the associated sets are both magog matrices and ASMs.  In contrast, (2) shows there is only one magog matrix with exactly one inversion, while there are $n-1$ ASMs with exactly one inversion. Finally, (5) gives a linear formula for magog matrices with inversion number two, whereas the number of ASMs with inversion number two grows faster.
\end{remark}
We give the following conjecture, based on the numerical data in Table~\ref{tab:MMinversion_stats}, computed in SageMath~\cite{sage}. Note that the second sequence matches OEIS~\cite{oeis} sequence A090809, the coefficient of a certain irreducible character of the symmetric group in a certain Kronecker power.
\begin{conjecture}
The following enumerations hold:
\begin{itemize}
    \item $|\magog_n(1~\posinv)|=\binom{n}{2}$, $n\geq 2$,
    \item $|\magog_n(2~\posinv)|=2\binom{n-1}{2} + 4\binom{n-1}{3} + 3\binom{n-1}{4}$, $n\geq 3$.
    \item $|\magog_n(\binom{n}{2}-2~\posinv)|=n(n-2)$, $n\geq 3$.
        \item $|\magog_n(\binom{n}{2}-1~\inv)|=2^n-n-1$, the $n$th Eulerian number.
\end{itemize}
\end{conjecture}

\begin{table}[hbtp]
\centering   
\begin{tabular}{|c | c | c | c |} 
 \hline
 $n$   & Distribution of   & Distribution of \\
      & positive inversions & inversions \\
\hline\hline
     3  & 1,3,2,1 & 1,1,4,1 \\
\hline
 4   & 1,6,10,13,8,3,1 & 1,1,5,11,12,11,1  \\
 \hline
 5  & 1, 10, 31, 70, 94, 90, & 1, 1, 6, 14, 39, 58, \\
      & 74, 39, 15, 4, 1 & 104, 105, 74, 26, 1\\
 \hline
 6 & 1, 15, 75, 259, 577, 954,  & 1, 1, 7, 17, 52, 132, \\
    & 1315, 1391, 1171, 829, 501, &  
   275, 541, 921, 1332, 1481,  \\
      &  234, 84, 24, 5, 1 &  
 1420, 856, 342, 57, 1 \\
 \hline
7 &1, 21, 155, 764, 2516, 6240, 12757,  & 1, 1, 8, 20, 66, 181, 509, 1139, 2573,\\
&21033, 28567, 33326, 33789,29256, 21730, 13983,  &  5275, 9970, 16752, 25117,33291, 37866, \\
 &  7909, 3935, 1619, 552, 153, 35, 6, 1&  35740, 26797, 15694, 5873, 1354, 120, 1\\
\hline
\end{tabular}
    \caption{Inversion data for magog matrices.}
    \label{tab:MMinversion_stats}
    \begin{tabular}{|c | c | c | } 
 \hline
 $n$   & Distribution of   & Distribution \\
      & positive inversions & of inversions \\
\hline\hline
     3  & 1,3,2,1 & 1,2,3,1  \\
\hline
 4    & 1,6,11,13,7,3,1 & 1,3,7,13,11,6,1  \\
 \hline
 5  & 1, 10, 35, 77, 99, 92, 67, 31, 12, 4, 1 & 1, 4, 12, 31, 67, 92, 99, 77, 35, 10, 1\\
 \hline
  6& 
  1, 15, 85, 302, 684, 1122, 1443,
 1396, 

  & 
  1, 5, 18, 56, 156, 379, 696, 1077, \\
  &1077, 696, 379, 156, 56, 18,
 5, 1
  &1396, 1443, 1122, 684, 302, 85, 15, 1\\
  \hline
  7
  &1, 21, 175, 917, 3196, 8166,
 16421, 
  & 1, 6, 25, 89, 287, 849, 2261, 5016, 9778, \\
  &26064, 33489, 36010, 33089, 25580,
 16907, 
 &16907, 25580, 33089,
 36010, 33489,\\
 & 9778, 5016, 2261, 849, 287,
 89, 25, 6, 1  
 & 26064, 16421, 8166, 3196,
 917, 175, 21, 1 \\
 \hline
\end{tabular}
    \caption{Inversion data for alternating sign matrices, for comparison to Table~\ref{tab:MMinversion_stats}.}
    \label{tab:ASMinversion_stats}
\end{table}

\section{TSSCPP matrix polytope}
\label{sec:magogpoly}
In this section, we define and study a polytope by taking the convex hull of $n\times n$ magog matrices. 

\subsection{TSSCPP matrix polytope definition and partial inequality description}
\label{subsec:magogpoly}
In this subsection, we use the matrices of Section~\ref{sec:matrices} to define a polytope, in an analogous manner as was done for alternating sign matrices in \cite{Striker2007TheAS}. In the next subsection, we give a partial inequality description of this polytope.

\begin{definition}
    We define the \emph{magog matrix polytope} of order $n$, denoted $\TSSCPP(n)$, to be the convex hull of all $n\times n$ magog matrices $\magog_n$.
\end{definition}

Our first  result on this polytope is an upper bound on the dimension.

\begin{proposition}
    The dimension of $\TSSCPP(n)$ is at most $(n-1)^2$.
\end{proposition}
\begin{proof}
    The ambient dimension of $\TSSCPP(n)$ is $n^2$, since there are $n^2$ entries in each matrix. But since the rows and columns sum to one, the last entry in each row and each column is determined by the rest. So the dimension is at most $(n-1)^2$. 
\end{proof}

\begin{remark}
    Our computational evidence indicates $(n-1)^2$ is the dimension, not merely an upper bound. Since only the $132$-avoiding permutation matrices are included in $\magog_n$, it is not true that $\TSSCPP(n)$
contains the Birkhoff polytope, which has dimension $(n - 1)^2$. This argument therefore cannot be used to show
$\TSSCPP(n)$ has dimension at least $(n - 1)^2$, as has been previously argued for the ASM polytope. We also cannot use the same proof method as for the Birkhoff polytope dimension, given in \cite{yemelichev1984polytopes}, since this relies on knowing all the inequalities. It may be sufficient to find the inequalities for the convex hull of the permutation matrices that are magog matrices, as the \emph{$132$-pattern avoiding polytope} appears to have dimension $(n-1)^2$; see \cite[Table 1]{SaganDavisPolytope}.
\end{remark}

\subsection{Inequalities and a partial facet description}
\label{subsec:magogpolyineq}
In the cases of the Birkhoff polytope and the ASM polytope~\cite{Behrend2007HigherSA,Striker2007TheAS}, the inequalities defining the matrices, when allowed to apply to real numbers rather than integers, also defined their convex hull. Since $\TSSCPP(n)$ is also a polytope defined as the convex hull of a subset of square sign matrices, one may expect the inequality description of $\TSSCPP(n)$ to be given by \eqref{eq:1magog}-\eqref{eq:4magog} and \eqref{eq:5magog}. This is not the case, as seen in the following remark.

\begin{remark}
\label{remark:ineq}
  Note that unlike in the case of permutation matrices or ASMs, the polytope defined using the inequalities satisfied by magog matrices, \eqref{eq:1magog}-\eqref{eq:4magog} and \eqref{eq:5magog}, is not $\TSSCPP(n)$. For example, when $n=3$ the $\mathcal{H}$-polytope defined by the inequalities given in \eqref{eq:1magog}-\eqref{eq:4magog} and \eqref{eq:5magog} has vertices including:
  \[\begin{array}{ccl}
      \begin{pmatrix} 0.5 & 0 & 0.5  \\ 0.5 & 0 & 0.5 \\ 0 & 1 & 0  \end{pmatrix}

    \hspace{1cm}
       \begin{pmatrix}  0.5 & 0.5 & 0  \\ 0 & 0 & 1 &  \\ 0.5& 0.5 & 0   \end{pmatrix}
\end{array}\]
    neither of which are magog matrices. 
    
    We have computed with SageMath that a complete inequality description of $\TSSCPP(3)$ is given by: $a_{1,1} \geq 0, a_{1,2} \geq 0, a_{1,3} \geq 0, a_{3,1} \geq 0, a_{3,2} \geq 0$, and $a_{1,2}+a_{2,1}+a_{2,2} \geq 1$. The last of these inequalities appears in \eqref{eq:7} of our partial inequality description in the next theorem.
\end{remark}

Though the naive guess of an inequality description for the  {magog matrix polytope} $\TSSCPP(n)$ does not work, in the theorem below, we give
 some additional inequalities satisfied by $\TSSCPP(n)$. Remark~\ref{remark:newineq} discusses why adding these to the naive guess also does not form a complete inequality description.

    \begin{theorem}
    \label{thm:magogineq}
       The following inequalities hold for any matrix $A=(a_{ij}) \in \TSSCPP(n)$: 
        \begin{align}
        \label{eq:7}        \sum_{j'=1}^{j}a_{i+1,j'}+\sum_{i'=1}^{i+1}a_{i',j+1} \geq 1
       &\text{ for all }
        1 \leq i \leq n-2,
        1 \leq j \leq n-2,
        i+j \geq n-1,  \\ 
        \label{eq:8} 
        \sum_{j'=1}^{j+1}a_{2,j'}+\sum_{j'=j+1}^{n}a_{1,j'} \geq 1  
        &\text{ for all } 
        1 \leq j \leq n-3, \\
        \label{eq:9}
\sum_{i'=1}^{i+1}a_{i',2}+\sum_{i'=i+1}^{n}a_{i',1} \geq 1        
        &\text{ for all }
        1 \leq i \leq n-3.    
         \end{align}
    \end{theorem}
    \begin{proof}
        It suffices to show the above inequalities are true for all magog matrices, thus they hold for convex combinations of magog matrices as well.
        
        For each set of inequalties, we show that if a $\{0,1,-1\}$-matrix $A$ violates that inequality, then $A$ is not a magog matrix. We start with \eqref{eq:7}. That is, we consider $i,j$ where         $1 \leq i \leq n-2$, 
        $1 \leq j \leq n-2$, 
        $i+j \geq n-1$,   and   
\begin{equation}
\sum_{j'=1}^{j}a_{i+1,j'}+\sum_{i'=1}^{i+1}a_{i',j+1} < 1.
\end{equation}
But in fact, since the left hand side is the sum of partial row and column sums of a magog matrix, it must be a nonnegative integer. Thus,
\begin{equation}
\sum_{j'=1}^{j}a_{i+1,j'}+\sum_{i'=1}^{i+1}a_{i',j+1} =0. \label{eq:lshape}
\end{equation}

 Note there are $i+1$ columns with positive partial sum $\sum_{i'=1}^{i+1}a_{i',j+1}$.
 Since $i+j \geq n-1$ so that $n-j\leq i+1$, there are also at least $n-j$ columns with positive partial sum $\sum_{i'=1}^{i+1}a_{i',j+1}$.

At least one of them must occur in the first $j+1$ columns so 
\[
            \sum_{i'=1}^{i+1}\sum_{j'=1}^{j+1}a_{i',j'} \geq 1.
\]        
We can subtract \eqref{eq:lshape} from both sides to obtain:
\[         
               \sum_{i'=1}^{i}\sum_{j'=1}^{j}a_{i',j'} \geq 1.
\]     
          This means that at least one column less than or equal to $j$ has a partial sum of one up to row $i$. We choose the largest column $\widetilde{j}$ with this property. Then
\[
         \sum_{i'=1}^{i}a_{i',\widetilde{j}} \geq 1.
\]
         
         If $\widetilde{j}=j$ we will violate the $(i,j)$-special inequality because of \eqref{eq:lshape}. Otherwise,
\[
         \sum_{i'=1}^{i}a_{i',\widetilde{j}+1} \leq 0     
\]       
         since otherwise we would have that $\widetilde{j}$ is not the last column before $j$ with a positive partial sum. Likewise, 
\[   \sum_{j'=1}^{\widetilde{j}+1}a_{i+1,j'} \leq 0
\]      
          otherwise we would have that there is a $-1$ between $\widetilde{j}+1$ and $j$ in row $i+1$ but that column would have a one earlier, meaning $\widetilde{j}$ is not the last partial sum of one. This means the $(i,\widetilde{j})$-special inequality is violated.

        Next we will show the same for \eqref{eq:8}. Suppose it does not hold. In the first row the one must occur between column one and $j$ in what we call column $\widetilde{j}$. For the same reasons discussed in the above paragraph the $(1,\widetilde{j})$-special inequality is violated.

        Next we prove \eqref{eq:9} holds. By subtracting $\sum_{i'=1}^{n}a_{i'1}=1$ from both sides, we can rewrite this as 
\[
        \sum_{i'=1}^{i+1}a_{i',2}-\sum_{i'=1}^{i}a_{i',1} \geq 0.
\]
        Now if this inequality does not hold, the partial sum of the first column up to $i$ must be one. This means that $a_{i+1,1}=0$ and the $(i,1)$-special inequality is violated.
    \end{proof}

\begin{remark}
\label{remark:newineq}
    If we were to define a polytope from the inequalities in \eqref{eq:1magog}-\eqref{eq:4magog}, \eqref{eq:5magog}, and \eqref{eq:7}-\eqref{eq:9}, we still do not have a complete inequality description of the TSSCPP polytope $\TSSCPP(n)$. In the case $n=4$, we have computed with SageMath that the inequality $a_{1,1} \leq a_{4,4}$ is also necessary.
\end{remark}

\subsection{Vertices}
\label{subsec:magogpolyvert}
In this subsection, we show the magog matrices are the extreme points of the polytope $\TSSCPP(n)$.
    
    \begin{theorem}
    \label{thm:magogextreme}
        The vertices of $\TSSCPP(n)$ are the set of $n \times n$ magog matrices. 
    \end{theorem}
    \begin{proof}
     Fix a magog matrix $A$. We show that there is a hyperplane separating $A$ from all other magog matrices. 
     Take $C_A$ as the set where $(i,j) \in C_A$ if the partial sum up to row $i$ in column $j$ is one. Notice that since there is one positive partial sum added for each row, every matrix has $\binom{n}{2}$ elements in $C_A$. 
     We consider the hyperplane $H_A(X)=\binom{n}{2}-\frac{1}{2}$ for $H(A)$ defined as follows:\[H_A(X)=\sum_{(i,j) \in C_A}\sum_{i'=1}^{i}x_{i',j}. \] 
     Now for $X=A$, $H_A(X)=\binom{n}{2}$. Given a hyperplane defined in this manner, we may recover the matrix from which it is formed, so $H_A$ is unique for each $A$. 
     By definition, every magog matrix has $\binom{n}{2}$ partial column sums that equal one. Now let $A'\neq A$ be another magog matrix; there must be $(i,j) \in C_A$ where where the column partial sum is one in $A$ and zero in $A'$. Then $H_A(A')$ will be less than $H_A(A')$ for every $i,j$ where this occurs. Therefore, $H_A(A)>\binom{n}{2}-1/2$ but $H_A(A')<\binom{n}{2}-1/2$.
     \end{proof}

We end this subsection with some remarks on data about $\TSSCPP(n)$ in Table \ref{tab:matrixpolytopestats1} as compared to the data for the ASM polytope in Table \ref{tab:matrixpolytopestats2}, both computed with SageMath.
\begin{remark}
Note from Tables \ref{tab:matrixpolytopestats1} and \ref{tab:matrixpolytopestats2} that for $n\leq 4$, the $f$-vector of the ASM polytope is  entrywise greater than or equal to that of the TSSCPP polytope. However for $n \geq 5$ the TSSCPP polytope has more facets and other high dimensional faces. Another interesting phenomena is that the diameter of the TSSCPP polytope equals two for all $2\leq n\leq 5$, but the diameter of the ASM polytope is $3$ when $n=5$.
\end{remark}
    
\begin{table}[hbtp]
\centering
\begin{tabular}{|c | c c c|} 
 \hline
 $n$ & Dimension & Normalized Volume  & $f$-vector   \\ 
 \hline\hline
 2 &  1 & 1 & 1,2,1 \\ 
 \hline
 3 & 4 & 3 & 1,7,15,14,6,1   \\
 \hline
 4 & 9 & 596 & 1, 42, 275, 809, 1355, 1421, 967, 425, 115, 17, 1
 \\
 \hline
 5 & 16 & 8558 &1, 429, 8558, 65181, 276977, 763114,  \\

&&&1471355, 2072494, 2185387,  1746405, 1060656, \\
&&& 486723, 166180, 41053, 7001, 760, 45, 1\\
\hline
\end{tabular}
\begin{tabular}{|c | c c  |} 
 \hline
 $n$   & Diameter & Ehrhart polynomial \\ 
 \hline\hline
 2  & 1 & $t+1$\\ 
 \hline
 3 & 2 & $(1/8)t^4 + (11/12)t^3 + (19/8)t^2 + (31/12)t + 1$  \\
 \hline
 4 & 2 & $(149/90720)t^9 + (443/13440)t^8 + (1237/4320)t^7 +$  
 \\
 &&$(1357/960)t^6 + (18899/4320)t^5 + (5647/640)t^4 + $\\
 &&$(1055983/90720)t^3 + (32693/3360)t^2 + (1691/360)t + 1$\\
 \hline
 5  & 2 & * \\

\hline
\end{tabular}
\caption{Data for the TSSCPP matrix polytope.}
    \label{tab:matrixpolytopestats1}
\end{table}
\begin{table}[hbtp]
\centering
\begin{tabular}{|c | c c c|} 
 \hline
 $n$ & Dimension & Normalized Volume  & $f$-vector   \\ 
 \hline\hline
 2 &  1 & 1 & 1,2,1 \\ 
 \hline
 3 & 4 & 4 & 1,7,17,18,8,1   \\
 \hline
 4 & 9 & 1376 & 1, 42, 380, 1333, 2428, 2580, 1675, 668, 158, 20, 1
 \\
 \hline
 5 & 16 & * &1, 429, 15841, 138583, 602644, 1610815,   \\

&&& 2920505, 3791870, 3638846, 2626814, 1436743, \\
&&& 594222, 183682, 41442, 6550, 676, 40, 1\\
\hline
\end{tabular}
\begin{tabular}{|c | c c  |} 
 \hline
 $n$   & Diameter & Ehrhart polynomial \\ 
 \hline\hline
 2  & 1 & $t+1$\\ 
 \hline
 3 & 2 & $(1/6)t^4 + t^3 + (7/3)t^2 + (5/2)t + 1$  \\
 \hline
 4 & 2 & $ (43/11340)t^9 + (311/5040)t^8 + (167/378)t^7 +  $  
 \\
 &&$(661/360)t^6 + (2659/540)t^5 + (6437/720)t^4 +$\\
 &&$(6301/567)t^3 + (962/105)t^2 + (1423/315)t + 1$\\
 \hline
 5  & 3 & * \\

\hline
\end{tabular}
\caption{Data for the ASM polytope, for comparison to Tables~\ref{tab:matrixpolytopestats1} and \ref{tab:booleanpolytopestats}.}
    \label{tab:matrixpolytopestats2}
\end{table}
\begin{remark}
 For the $n$th ASM polytope, the number of facets is $4n^2-16n+20$ for $n \geq 3$ \cite[Theorem 3.3]{Striker2007TheAS}. Now if we have a polynomial of degree $k$ and we know $k+1$ points on the polynomial, those facts define the polynomial uniquely. Thus if we suppose that the TSSCPP matrix polytope facet formula $f(x)$ is a quadratic for all $n \geq 3$ we know $f(3)=6$, $f(4)=17$, and $f(5)=45$. The quadratic defined by these values is $f(x)=8.5x^2-48.5x+75$. However, $f(6)=90$, but the six dimensional TSSCPP matrix polytope has $162$ facets, as confirmed by SageMath computations. This suggests that if the TSSCPP matrix polytope facet counting formula is a polynomial, it has order higher than two.
\end{remark}

\section{TSSCPP boolean triangle polytope}
\label{sec:btp}
In this section, we study a polytope defined as the convex hull of different objects in bijection with TSSCPPs. Recall Definition~\ref{def:bool} and Examples~\ref{ex:bool1} and \ref{ex:bool2}.

\begin{definition}   
 Define the TSSCPP boolean triangle polytope $\BTP(n)$ as the subset of $\RR^{\binom{n}{2}}$ that is the convex hull of all TSSCPP boolean triangles of order $n$.
\end{definition}

\begin{proposition}
    The TSSCPP boolean triangle polytope $\BTP(n)$ is of dimension $\binom{n}{2}$.
\end{proposition}
\begin{proof}
  The dimension of this polytope is $\binom{n}{2}$ since all the  standard basis vectors of $\RR^{\binom{n}{2}}$ as well as the zero vector are TSSCPP boolean triangles.
\end{proof}

We use a proof inspired by the Birkhoff-von Neumann theorem to obtain the following inequality description of $\BTP(n)$. See Example~\ref{ex:boolineqex} for an example of the proof method.
    \begin{theorem}
    \label{thm:btpineq}
        The TSSCPP boolean triangle polytope $\BTP(n)$ is equal to the set of all triangular arrays $(b_{i,n-j})$ for $1 \leq j \leq i \leq n-1$ such that: 
        \begin{align}
        b_{i,n-j} &\geq 0, \label{eq:11}\\
        b_{i,n-j} &\leq 1, \label{eq:12}\\
    1+\sum_{k=j+1}^{i}b_{k,n-j-1} &\geq \sum_{k=j}^{i}b_{k,n-j}   \label{eq:13} \text{ for all } 1 \leq j < i \leq n-1. 
    \end{align}
    We refer to \eqref{eq:13} as the $(i,j)$-diagonal inequality.
    \end{theorem}

    \begin{proof}
        Call the set defined by the above inequalities $P(n)$. It is easy to see that any convex combination of TSSCPP boolean triangles is in $P(n)$. Now we show any $B=(b_{ij}) \in P(n)$ is a convex combination of boolean triangles. Our goal is to show that $B$ can be written as a convex combination of two members of $P(n)$ with at least one fewer non-integer entry or one fewer non-integer partial sum difference. From there induction proves the rest and we are done.

        We label some elements of our triangle in the following manner: if $\sum_{k=j}^{i-1}b_{k,n-j} \in \ZZ$ but $\sum_{k=j}^{i}b_{k,n-j} \notin \ZZ$ (meaning $b_{ij}$ takes the diagonal partial sum from the integers to outside the integers) we label $b_{ij}$ with a $(+)$. On the other hand if $\sum_{k=j}^{i-1}b_{k,n-j} \notin \ZZ$ but $\sum_{k=j}^{i}b_{k,n-j} \in \ZZ$ then we label $b_{ij}$ with a $(-)$. We consider a northwest to southeast diagonal at an entry $b_{ij}$ to be in an \emph{unbalanced} state if there are
an unequal number of entries labelled $(+)$ and $(-)$ up to and including $b_{ij}$, and in a \emph{balanced} state
when the numbers are equal.

        We define $k'=\min\{1-x_{ij},y_{ij},1-p_{ij}\}$  
        where the $x_{ij}$ denote all of the entries labeled $(+)$, the $y_{ij}$ denote all of the entries labeled $(-)$, and the $p_{ij}$ are all of the differences in diagonal partial sums $\left(\sum_{k=j}^{i}b_{k,n-j}-\sum_{k=j+1}^{i}b_{k,n-j-1}\right)$ when $b_{ij}$ is in a balanced state and $b_{i,j+1}$ is in an unbalanced state. Now we subtract $k'$ from from all entries labeled $(-)$ and add $k'$ to all entries labeled $(+)$ to get $B'$. This keeps all entries between zero and one, keeps all partial sum differences less than or equal to one, and adds one more integer entry or partial sum difference. 
        
        We show that this triangle is still in $P(n)$. For each of the $(i,j)$-inequalities of \eqref{eq:13}, there are four cases: the $j$ diagonal and $j-1$ diagonal are both in an unbalanced state, both diagonals are in a balanced state, the $j$ diagonal is in a balanced state and $j-1$ is in the unbalanced state, or the $j$ diagonal is in an unbalanced state and the $j-1$ diagonal is in a balanced state. In the first two cases, the difference in partial sums stays constant, in the third case the partial sum difference can only decrease, and in the fourth case   $k'$ is bounded by $1-p_{ij}$. All entries will remain between zero and one since we bound $k'$ by $1-x_{ij}$ and $y_{ij}$. 
        
        Now we switch the $(+)$ and $(-)$ and define $k''=\min\{1-x_{ij},y_{ij},1-q_{ij}\}$, where the $q_{ij}$ are all of the differences in diagonal partial sums 
        when $b_{ij}$ is in an unbalanced state and $b_{i,j+1}$ is in a balanced state. Subtract $k''$ from $(-)$ and add to $(+)$ to get $B''$, which we may show in an analogous way as for $B'$ is also in $P(n)$. Then $B=\frac{k''}{k'+k''}B'+\frac{k'}{k'+k''}B''$. 
    \end{proof}
    
    \begin{example}
    \label{ex:boolineqex}
The first triangle below is an example element of $P(n)$ on which we apply the algorithm of Theorem~\ref{thm:btpineq} to create a convex combination. The second triangle demonstrates which elements we are changing to get our convex combination. The blue elements are the northwest to southeast diagonal partial sums. Each red number is the difference between the blue diagonal partial sums to its right and left. The circled elements are in an unbalanced state. The squared elements are the elements that are in a balanced state. The third diagram demonstrates the resulting two triangles that give  our convex combination. To see this in practice, the largest element whose partial sum is a $(+)$ is $0.8$, the smallest with a $(-)$ is $0.2$ and largest difference in partial sums between unbalanced and balanced states is $0.5$. This means that $k'=\min\{1-0.8,0.2,1-0.5\}=0.2$. Likewise, the largest element whose partial sum is a $(-)$ is $0.9$ and the smallest whose is a $(+)$ is 0.1, the largest difference in partial sums is 0.9. Thus $k''=\min\{1-0.9,0.1,1-0.9\}=0.1$ and we obtain the desired result.
   \begin{center}
        $\begin{array}{cccccccccc}
         & & &  &  &   0.5 &  &    & & \\
         & &  &  & 0.8 &  & 0 &  &  &  \\ 
         & & & 0.1 &  & 0.2 &  & 1 &  &   \\
         & & 1 &  & 0.9 &  & 1 &  & 0.5 &  \\
         & 0.1 & &  0.1 & & 0.1 & & 0.1& &1 \\    
     \end{array}$~\hspace{-2.5in}\raisebox{-1.2in}{\includegraphics{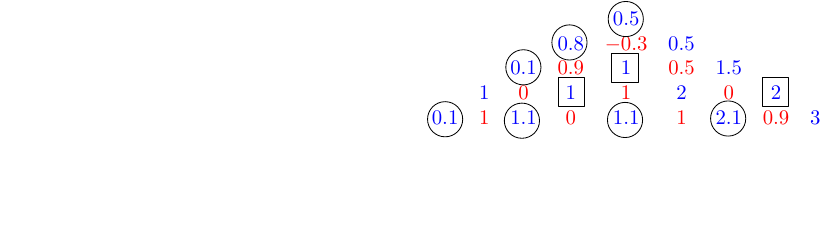}}

\vspace{-.5in}
  \scalebox{.70}{$\begin{bmatrix}
             & & &  &  &   0.5 &  &    & & \\
         & &  &  & 0.8 &  & 0 &  &  &  \\ 
         & & & 0.1 &  & 0.2 &  & 1 &  &   \\
         & & 1 &  & 0.9 &  & 1 &  & 0.5 &  \\
         & 0.1 & &  0.1 & & 0.1 & & 0.1& &1 \\
        \end{bmatrix}
        =
        \displaystyle\frac{0.1}{0.3}
        \begin{bmatrix}
             & & &  &  &   0.7 &  &    & & \\
         & &  &  & 1 &  & 0 &  &  &  \\ 
         & & & 0.3 &  & 0 &  & 1 &  &   \\
         & & 1 &  & 0.7 &  & 1 &  & 0.3 &  \\
         & 0.3 & &  0.3 & & 0.3 & & 0.3& &1 \\
        \end{bmatrix}
        +
        \displaystyle\frac{0.2}{0.3}
        \begin{bmatrix}
             & & &  &  &   0.4 &  &    & & \\
         & &  &  & 0.7 &  & 0 &  &  &  \\ 
         & & & 0 &  & 0.3 &  & 1 &  &   \\
         & & 1 &  & 1 &  & 1 &  & 0.6 &  \\
         & 0 & &  0 & & 0 & & 0& &1 \\
        \end{bmatrix}$
 }
 \end{center}
    \end{example}

Next we use this inequality description to enumerate the facets of the TSSCPP boolean triangle polytope.   
    \begin{theorem}
    \label{thm:btpfacets}
        The number of facets of $\BTP(n)$ is $\displaystyle\frac{(n-1)(3n-2)}{2}$.
    \end{theorem}
    \begin{proof}
        In the previous theorem we give an inequality description with $2\binom{n}{2}+\binom{n-1}{2}=\frac{(n-1)(3n-2)}{2}$ inequalities; this gives an upper bound for the number of facets. To show that this is also a lower bound, for each inequality, we find a vertex that lies on that inequality but no other.

        First we demonstrate this for inequalities of the form \eqref{eq:11}. Fix $i',j'$ and consider the triangular array $\{a_{ij}\}$ defined by the equation below (an example is drawn in triangular form to the right for $n=6$, $i'=3$, $j'=4$):
        \begin{multicols}{2}
        \begin{equation*}
        a_{ij}=
        \begin{cases}
        0& i=i' \text{ and } j=j'\\
        0.25 & j=j'+1 \text{ and } i=i' \\
        0.5 & \text{otherwise.}
        \end{cases}    
        \end{equation*}   
                $\begin{array}{ccccccccc}
          & &  &  &   0.5 &  &    & & \\
          &  &  & 0.5 &  & 0.5 &  &  &  \\ 
          & & 0.5 &  & 0 &  & 0.25 &  &   \\
          & 0.5 &  & 0.5 &  & 0.5 &  & 0.5 &  \\
          0.5 & &  0.5 & & 0.5 & & 0.5& &0.5   
     \end{array}$
        \end{multicols}
It is easy to see that \eqref{eq:11} and \eqref{eq:12} do not hold with equality for any other $i,j$.  The difference between the partial sums will be $0.5$ when comparing $j'$ and $j'-1$ before we reach $i'$ and zero after we reach $i'$. For $j'$ and $j'+1$ the difference will be $0.5$ before we reach $i'$ and $0.75$ after we reach $i'$. For $j'+1$ and $j'+2$ the difference is $0.5$ before $i'$ and $0.75$ after $i'$. For all other diagonals the difference is $0.5$. 

 Next we  demonstrate this for inequalities of the form \eqref{eq:12}. Fix $i',j'$ and consider the triangular array $\{a_{ij}\}$ defined by the equation below (an example is drawn in triangular form to the right for $n=6$, $i'=3$, $j'=4$):
        \begin{multicols}{2}
        \begin{equation*}
a_{ij}=
\begin{cases}
    1& i=i' \text{ and } j=j'\\
      0.75 & i=i' \text{ and } j=j'-1 \\
      0.5 & \text{ otherwise.}
   \end{cases}
\end{equation*}
$\begin{array}{ccccccccc}
          & &  &  &   0.5 &  &    & & \\
          &  &  & 0.5 &  & 0.5 &  &  &  \\ 
          & & 0.75 &  & 1 &  & 0.5 &  &   \\
          & 0.5 &  & 0.5 &  & 0.5 &  & 0.5 &  \\
          0.5 & &  0.5 & & 0.5 & & 0.5& &0.5   
     \end{array}$
 \end{multicols}
Agagin, it is easy to see that \eqref{eq:11} and \eqref{eq:12} do not hold with equality for any other $i',j'$. The difference between the partial sums will be $0.5$ when comparing $j'$ and $j'-1$ before we reach $i'$ and $0.75$ after we reach $i'$. For $j'$ and $j'+1$ the difference will be $0.5$ before we reach $i'$ and zero after we reach $i'$. For $j'-1$ and $j'-2$ the difference is $0.50$ before we reach $i'$ and $0.75$ after we reach $i'$. For all other diagonals the difference is $0.5$. 

Finally we  demonstrate this for inequalities of the form \eqref{eq:13}. Fix $i',j'$ and consider the triangular array ${a_{ij}}$ defined by the equation below (an example is drawn in triangular form to the right for $n=6$, $i'=4$, $j'=4$):
\begin{multicols}{2}
        \begin{equation*}
a_{ij}=
\begin{cases}
         0.75& i=n-j' \text{ and } j=j'\\
      0.75 & i=i' \text{ and } j=j' \\
       0.25 & i=i'+1 \text{ and } j=j' \\
      0.5 & \text{otherwise}
   \end{cases}
\end{equation*} 
 \hspace{-1ex}$\begin{array}{ccccccccc}
          & &  &  &   0.5 &  &    & & \\
          &  &  & 0.75 &  & 0.5 &  &  &  \\ 
          & & 0.5 &  & 0.5 &  & 0.5 &  &   \\
          & 0.5 &  & 0.5 &  & 0.75 &  & 0.5 &  \\
          0.5 & &  0.5 & & 0.5 & & 0.25& &0.5   
     \end{array}$
\end{multicols}
It is easy to see that \eqref{eq:11} and \eqref{eq:12} do not hold for any $i',j'$. The difference between the partial sums will be $0.75$ when comparing $j'$ and $j'-1$ before we reach $i'$, one after we reach $i'$, and $0.75$ after $i'+1$. For $j'$ and $j'+1$ the difference will be $0.75$ before we reach $i'$, $0.5$ after we reach $i'$, and $0.75$ after we reach $i'+1$. For all other diagonals the difference is $0.5$. Thus, \eqref{eq:13} only holds when $i=i'$ and $j=j'$. 
    \end{proof}

    It is interesting to note that this facet enumeration of $\frac{(n-1)(3n-2)}{2}$ is the formula for the second pentagonal numbers, sequence A005449 in the OEIS \cite{oeis}.

Finally, we show that each TSSCPP boolean triangle is an extreme point of $\BTP(n)$.
    
    \begin{theorem}
    \label{thm:btpextreme}
        The vertices of $\BTP(n)$ are the $\TSSCPP$ boolean triangles of order $n$.
    \end{theorem}
    \begin{proof}
        To prove this, we fix a boolean triangle $B=(b_{ij})$ of order $n$ and construct a hyperplane that has $B$ on one side and all other boolean triangles on the other side. From there it will follow that since we are taking the convex hull, the boolean triangles are the vertices.

        Let $A_B$ be the set of all $(i,j)$ such that $b_{ij}=1$. $A_B$ is unique for each $B$. Choose a boolean triangle $X=(x_{ij})$. Then we consider the hyperplane:
        \begin{equation}
            H_B(X)=\sum_{(i,j) \in A_B}x_{ij}-\sum_{(i,j) \notin A_B}x_{ij}.
        \end{equation}
        Now we know $H_B(B)=|A_B|$. For any other triangle $X$,  we show $H_B(X) < |A_B|$. If $H_B(X) \geq |A_B|$ then $X$ must have ones in at least all of the indices that $B$ does. In addition, $X$ may not have any ones outside of those that $B$ has.  Since any deviation in $X$ from $B$ will decrease $H_B(X)$ by at least one, $H_B(X)=|A_B|-0.5$ is the necessary hyperplane.
    \end{proof}

See Table~\ref{tab:booleanpolytopestats} for more combinatorial data on these polytopes.

\begin{table}
\begin{tabular}{|c | c c c |} 
 \hline
 $n$ & Dimension & Normalized Volume  & $f$-vector  \\ 
 \hline\hline
 2 & 1 & 1 & $1,2,1$  \\ 
 \hline
 3 & 3 & 5 & $1,7,12,7,1$  \\ 
 \hline
 4 & 6 & 410 & $1, 42, 171, 275, 214, 83, 15, 1$ 
  \\
 \hline
 5 & 10 & * & $1, 429, 3594, 11801, 20755, 21986, $ 
 \\
 &&&$14739, 6338, 1725, 285, 26, 1$\\
 \hline
\end{tabular}
\begin{tabular}{|c | c c|} 
 \hline
 $n$  & Diameter & Ehrhart polynomial \\ 
 \hline\hline
 2&1 &$t + 1$  \\ \hline
 3&2 &$(5/6)t^3 + (5/2)t^2 + (8/3)t + 1$  \\ 
 \hline
  4& 4 & $(41/72)t^6 + (41/12)t^5 + (319/36)t^4 + (38/3)t^3 + (761/72)t^2 + (59/12)t + 1$ \\
\hline
5&5& $(4496/14175)t^{10} + (8992/2835)t^9 + (884033/60480)t^8 + (61709/1512)t^7 +$\\
&& $(1651621/21600)t^6 + (217183/2160)t^5 + (17088191/181440)t^4 + $
\\
&& $(564799/9072)t^3 + (703609/25200)t^2 + (1943/252)t + 1$
 \\
 \hline
\end{tabular}
\caption{Data for the TSSCPP Boolean triangle polytope.}
    \label{tab:booleanpolytopestats}
\end{table}

\section*{Conflict of interest statement}
On behalf of all authors, the corresponding author states that there is no conflict of interest.

\section*{Acknowledgements}
The authors thank the anonymous referees for their detailed reading and helpful comments. We also thank Greta Panova for suggesting the Boolean triangle polytope and C\u{a}t\u{a}lin Ciuperc\u{a} for helpful comments on the exposition. We also thank the developers of SageMath and CoCalc, which were helpful in this research. Striker was supported by a Simons Foundation/SFARI gift (MP-TSM-00002802) and NSF grant DMS-2247089.

\bibliographystyle{plain}
\bibliography{biblio}

\end{document}